\documentclass[11pt]{article}
\usepackage{amsmath,amsthm,amssymb,amscd,graphicx,psfrag,subfigure}

\title{\bf Capture zones of the family of functions  $\lambda z^m \exp(z)$. }

\author{N\'uria Fagella \thanks{Supported by DGICYT Grant No PB96-1153, BFM2000-0805-C02-01 and CIRIT 2001SGR-70}\\
{\small Dep.~de Matem\`atica Aplicada i An\`alisi}\\
{\small Universitat de Barcelona}\\
{\small Gran Via de les Corts Catalanes, 585}\\
{\small 08005 Barcelona, Spain}\\
{\small \tt   fagella@maia.ub.es} \and Antonio Garijo \thanks{Supported by DGI Grant No AYA2001-0762 }\\
{\small Dep.~d'Eng. Inform\`atica i Matem\`atiques}\\
{\small Universitat Rovira i Virgili}\\
{\small Av. Pa\"{\i}sos Catalans, 26}\\
{\small 43007 Tarragona, Spain}\\
{\small \tt  agarijo@etse.urv.es}}

\theoremstyle{plain}
\newtheorem{theorem}{Theorem}[section]
\newtheorem{lemma}[theorem]{Lemma}
\newtheorem{proposition}[theorem]{Proposition}

\newtheorem*{thmA}{Theorem A}
\newtheorem*{thmB}{Theorem B}
\newtheorem*{thmC}{Theorem C}

\theoremstyle{definition}
\newtheorem*{definition}{Definition}

\newtheorem{remark}[theorem]{Remark}

\newsavebox{\savepar}

\newcommand{\bc}{{\mathbb C}}
\newcommand{\br}{{\mathbb R}}
\newcommand{\bz}{{\mathbb Z}}

\newcommand{\bd}{{\mathbb D}}
\newcommand{\bp}{{\mathbb T}}

\newcommand{\la}{\lambda}

\date{}

	{\end{list}}

\begin{document}

\maketitle

 \vspace{0.2cm}
 \centerline{\bf Abstract}
 \begin{quotation}
 {\noindent \small We consider the family of entire transcendental maps
given by $F_{\lambda,m} (z ) \, = \, \lambda z^m \exp(z) $  where $ m \ge 2$.
All functions $F_{\lambda,m}$  have a superattracting
fixed point at $z=0$, and a critical point at $z=-m$. In the dynamical plane
we study the topology of the basin of attraction of $z=0$. In the parameter plane we focus
on the capture behaviour, i.e., $\lambda$ values such that the critical point
belongs to the basin of attraction of $z=0$. In particular, we find a capture zone for
which this basin has a unique connected component, whose boundary is then non-locally connected. 
However, there are parameter values for which the boundary of the immediate basin of $z=0$ is a 
quasicircle.}
\end{quotation}


\section{Introduction}

Our goal in this paper is to study some dynamical aspects of the families of
entire transcendental maps
\[
F_{\lambda,m} (z ) \, = \, \lambda
z^m \exp(z), \quad  m \ge 2.
\]

Observe that $m=0$ corresponds to the exponential family $E_{\lambda} (z) \,
= \, \lambda \exp(z)$, the simplest example of an entire transcendental
map with a unique asymptotic value, $z=0$, in analogy with the well known quadratic family of
polynomials $z \to z^2 +c$. The exponential map has been thoroughly studied by many authors (see
for example, \cite{DK}, \cite{DT}).

The case $m=1$ corresponds to  $G_{\lambda}(z) \, = \, \lambda z \exp(z)$ which appeared for the first time in \cite{Ba1} as an example of an entire transcendental map whose Julia set is the whole plane
(for an appropriate value of $\la$). Later on this family was studied in \cite{F} and \cite{G}. The asymptotic value $z=0$ of $G_{\lambda}(z)$ is fixed and its multiplier is $G'_{\lambda}(0) \, = \, \lambda$. Hence its dynamical character depends on the parameter $\la$. Besides this point, the
dynamical behavior of  $G_{\lambda}$ is determined by the orbit of the critical point $z=-1$.

Some functions in the family $F_{\lambda,m}=\lambda z^m \exp(z)$ for $m \ge 2$ have been used in the literature
as examples of certain dynamical phenomena (see for example \cite{B}, for a Baker domain
at a positive distance from any singular orbit). But, to our knowledge, no systematic study
has been made before this work. 

All functions of the form $F_{\lambda,m}$, with $m \ge 2$,  have a superattracting fixed point 
at $z=0$ of multiplicity $m$, which is also an asymptotic value. The only other 
critical point is $z=-m$. The coexistence of a superattracting fixed point and a free critical
point makes this family very much analogous to the family of cubic polinomials  $C_a(z) \, = 
\, z^3 \, - \, 3a^2z + 2a^3+a $ which is described in \cite{M}.

Let $f$ be a transcendental entire function. It is known that the dynamical plane can be decomposed
into two invariant sets. The first one is an open set, namely the Fatou set or stable set, denoted by
 ${\cal F}(f)$, and it is formed by points, $z_0 \in \bc$ whose iterates
$\{f^{\circ n}\}$ form a normal family, in the sense of Montel, in some neihbourhood of $z_0$. The 
second one, namely the Julia or chaotic set, denoted by  ${\cal J} (f)$, is defined as the 
complement of the Fatou set, that is ${\cal J} (f) = \bc - {\cal F}(f)$. Fatou (\cite{Fa}) showed
that the Julia set of an entire transcendental function is a completely invariant, closed, nonempty 
and perfect set. As in the polynomial case, it may also be defined as the closure of the set of repelling periodic 
points.

We denote by $A(0)\,=\,A_{\lambda,m}(0)$ the basin of attraction of $z=0$, i.e., the set of all $z$ such that $F_{\lambda,m}^{\circ n}(z) \to 0 $ as $n$ tends to $+ \infty$. We also denote by $A^*(0)\,=\,A_{\lambda,m}^*(0)$ the 
connected component of $A(0)$ containing $z=0$. 

\begin{figure}[hbt]
\psfrag{d}[][]{\small $D_{\epsilon_0}$}
\psfrag{s0}[][]{\small $S^0_{\lambda,m}$}
\psfrag{s1}[][]{\small $S^1_{\lambda,m}$}
\psfrag{s2}[][]{\small $S^2_{\lambda,m}$}
\psfrag{sm1}[][]{\small $S^{-1}_{\lambda,m}$}
\psfrag{sm2}[][]{\small $S^{-2}_{\lambda,m}$}
\psfrag{H}[][]{\small $H_{|\lambda|,m}$}
\centerline{\includegraphics[width=0.4 \textwidth]{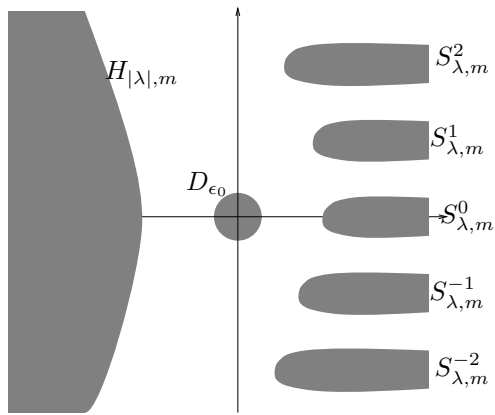}}
\caption{\small{Sketch of the basin of attraction of $z=0$, satisfying  that $D_{\epsilon_0} \subset
A^*(0)\, , \, F_{\lambda,m} (H_{|\lambda|,m}) \subset D_{\epsilon_0} \, \hbox{and} \,
 F_{\lambda,m} \, (S_{\lambda,m}^k) \subset H_{|\lambda|,m}$.}}
\label{a0}
\end{figure}

In Sec. 2  we concentrate on the dynamical plane and study the basin of attraction 
$A(0)$ (see Fig. \ref{a0}). The skeleton of the main components of $A(0)$  is needed to study later the parameter 
planes. We summarize the main results with regard to $A(0)$ in the following theorem.

\begin{thmA}
$\,$ \par
\begin{enumerate}
\item There exists $\epsilon_0 = \epsilon_0 (|\lambda|,m) \, > \, 0$, defined as the unique positive solution of $x^{m-1
} e^x = 1/|\lambda| $, such that $A^*(0)$ contains the 
disk $D_{\epsilon_{0}} \, = \, \{ z \in {\mathbb C} \, ; \, |z| \, < \, \epsilon_0 \} $. 
\item  There exist $x_0\,=\,x_0(|\lambda|,m) <0$ and a function $\, C(x) \ge 0$ such that the open set
 $$ H_{ |\lambda|,m} \, = \, \left\{ z \, = 
\, x \, + \, y i \, \left| \begin{array}{ll}
x \,&  \in \, (-\infty,x_0)  \\
y \,&  \in  (-C(x),C(x))  \end{array}\right. \right\}  $$ 

satisfies  $F_{\lambda,m} \, (H_{|\lambda|,m}) \, \subset \, D_{\epsilon_0}$.
\item There exist infinitely many strips, denoted  by $S_{\lambda,m}^k$, which are preimages of $H_{ |\lambda|,m}$. These
horizontal strips extend to $+\infty$, and they have asymptotic width equal to $\pi$.
\end{enumerate}
\end{thmA}

Section 3 is dedicated to parameter planes.  For some parameter values, the free critical point
$z=-m$ belongs to the basin of attraction of $z=0$, A(0), in which case we say that it has
been ``captured''. The connected components of parameter space for which this phenomenon occurs are
called {\it capture zones}, and they clearly do not exist for members of the family with $m<2$ such as
the exponential. Hence it is natural to study the set of parameters for which  the orbit of
 $z=-m$ is bounded. That is, we define the sets

$$B_m \, = \, \{ \lambda \in  {\mathbb C} \, | \, F_{\lambda,m}^{\circ n} (-m)  \nrightarrow
\infty \}.$$

In each of these sets, we may also distinguish between two different behaviours. Those parameter
values for which $-m \in A(0)$ and those for which this does not occur. Let $\stackrel{\circ}{B_m}$ denote the
interior of $B_m$. We will study the sets
\[
C_m^n = \{ \lambda \in \stackrel{\circ}{B_m} | F_{\lambda,m}^n (-m) \in A^*(0) \hbox{ and $n$ is the smallest number with this property} \}
\]

Although each $B_m$ contains infinitely many different capture zones, there is one which is
dynamically very different from all others.
We define the {\it main capture zone} $C_m^0$ as the set of parameter values $\lambda$ for which the
critical point $\,z=-m\,$ belongs to the immediate basin of $0$. That is,
\[
C_m^0 = \,  \{ \lambda \in \stackrel{\circ}{B_m} \, | \,  -m  \in A^*(0) \}.
\]

\begin{figure}[hbt]
\psfrag{ti1}[][]{\small $m \,$ even}
\psfrag{ti2}[][]{\small $m \,$ odd}
\psfrag{H}[][]{\small $C_m^2$}
\psfrag{D}[][]{\small $C_m^0$}
\psfrag{S}[][]{\small $C_m^3$}
\centerline{\includegraphics[width=0.8 \textwidth]{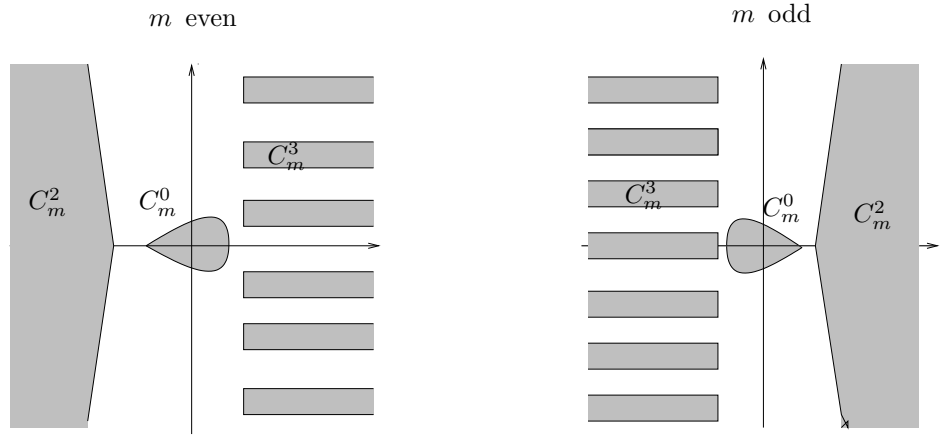}}
\caption{\small{Sketch of capture zones contained in $\stackrel{\circ}{B_m}$. On the left hand side when $m$ is even and on the right hand side when $m$ is odd.}}
\label{zona}
\end{figure}

As we shall see, this is a quite special component of $B_m$ since its boundary separates the parameter
values for which $\partial A^*(0)$ is a Cantor bouquet from those for which it is a Jordan curve (also,
this boundary separates the parameter values for which ${\cal F}(F_{\lambda,m})$ has one or infinitely
many components). The detailed study of this boundary will be the object of a future paper. Our goal
in Sec. 3  is to describe the main features of the parameter planes of the functions
$F_{\lambda,m}$ and, in particular, the structure of the capture zones (Fig. \ref{zona}). We 
summarize some of these facts in the following theorem.

\begin{thmB} 
$\,$ \par
\begin{enumerate}
\item The critical point $-m$ belongs to $A^*(0)$ if and only if the critical value $F_{\lambda,m}(-m)$ belongs 
to $A^*(0)$. Hence $C_m^1 \, = \, \o$.
\item The main capture zone $C_m^0$ is bounded.
\item The set $C_m^0$ contains the disk $\{ \lambda \in {\mathbb C} \, ;\, |\lambda| \, < \,
min(\frac{1}{e},(\frac{e}{m})^m ) \}$.
\item If $\lambda \in C_m^0$ then $A(0)=A^*(0)$, i.e., the basin of attraction of $z=0$ has a unique connected
component and hence it is totally invariant. Moreover, the boundary of $A^*(0)$ (which equals the Julia set) is a 
Cantor bouquet and hence it is disconnected and non-locally connected.
\item If $\lambda \notin C_m^0$ then $A(0)$ has infinitely many components. Moreover, if $|\lambda| > (\frac{e}{m-1})^{m-1} $,  the boundary of $A^*(0)$ is a quasi-circle.
\end{enumerate}
\end{thmB}

We also summarize some properties of the most obvious capture zones $C_m^2$ and $C_m^3$.

\begin{thmC}
$\,$ \par 
\begin{enumerate}
\item The set $C_m^2 $  contains an unbounded set to the left or to the right depending on the oddity of $m$. More
precisely, there exists a real constant $D_0(m) > 0$, and a function  $\alpha=\alpha(|\lambda|,m) \in (\pi/2,\pi)$, such that
\[
\begin{array}{ll}
\bullet \hbox{ for $m$ even,  the set $C_m^2$ contains the open set} &  \left\{ \lambda \in {\mathbb C} \, \left|
\begin{array}{ll}
 |\lambda| & \, >  \,  D_0   \\
|Arg(\lambda)|  & \, > \, \alpha    \end{array}\right. \right\} \\
\bullet \hbox{ for $m$ odd,  the set $C_m^2$ contains the open set} & \left\{ \lambda \in {\mathbb C} \, \left|
\begin{array}{ll}
|\lambda|  & \, > \, D_0   \\
 |Arg(\lambda)| & \, < \, \pi- \alpha  \end{array}\right. \right\}
\end{array}
\]
\item There exists infinitely many strips in $C_m^3$. If $m$ is even (resp. odd) then these horizontal strips
extend to $+\infty$  (resp. $-\infty$) and they have an asymptotic width equal to $(\frac{e}{m})^m \pi$.
\end{enumerate}
\end{thmC}

\section{The Dynamical Planes}

Our goal in this Sec. is to describe the dynamical plane of the families of maps
$F_{\lambda,m}$ given by the equation $F_{\lambda,m} (z ) \, = \, \lambda
z^m \exp(z), \quad \hbox{where } m \ge 2$.

The function $F_{\lambda,m}$ is a critically finite entire function, that is, it has a finite number of
 asymptotic values ($z=0$), and critical values ($z=0$ and $z=(-1)^m \lambda (\frac{m}{e})^m$). For this kind of
functions there exists a characterization of the Julia set (\cite{DT}), namely as the closure of the
set of points whose orbits tend to $\infty$.

Using the characterization above we can plot an
approximation of ${\cal J}(F_{\lambda,m})$. Generally, orbits tend to $\infty$ in specific directions. In
our case, if  $\lim_{n \to \infty} |F^{\circ n}_{\lambda,m} (z) | \, = \, +\infty \,$, then we have
$ \, \lim_{n \to \infty} Re(F^{\circ n}_{\lambda,m} (z) )  \, = \, +\infty$. Thus, an approximation
of the Julia set is given by the set of points whose orbit containts a point with real part greater
than, say, $50$.

\begin{figure}[hbt]
	\centering
	\subfigure[\scriptsize{$\lambda = -2.1 $} ]{
	\includegraphics[width=0.3 \textwidth]{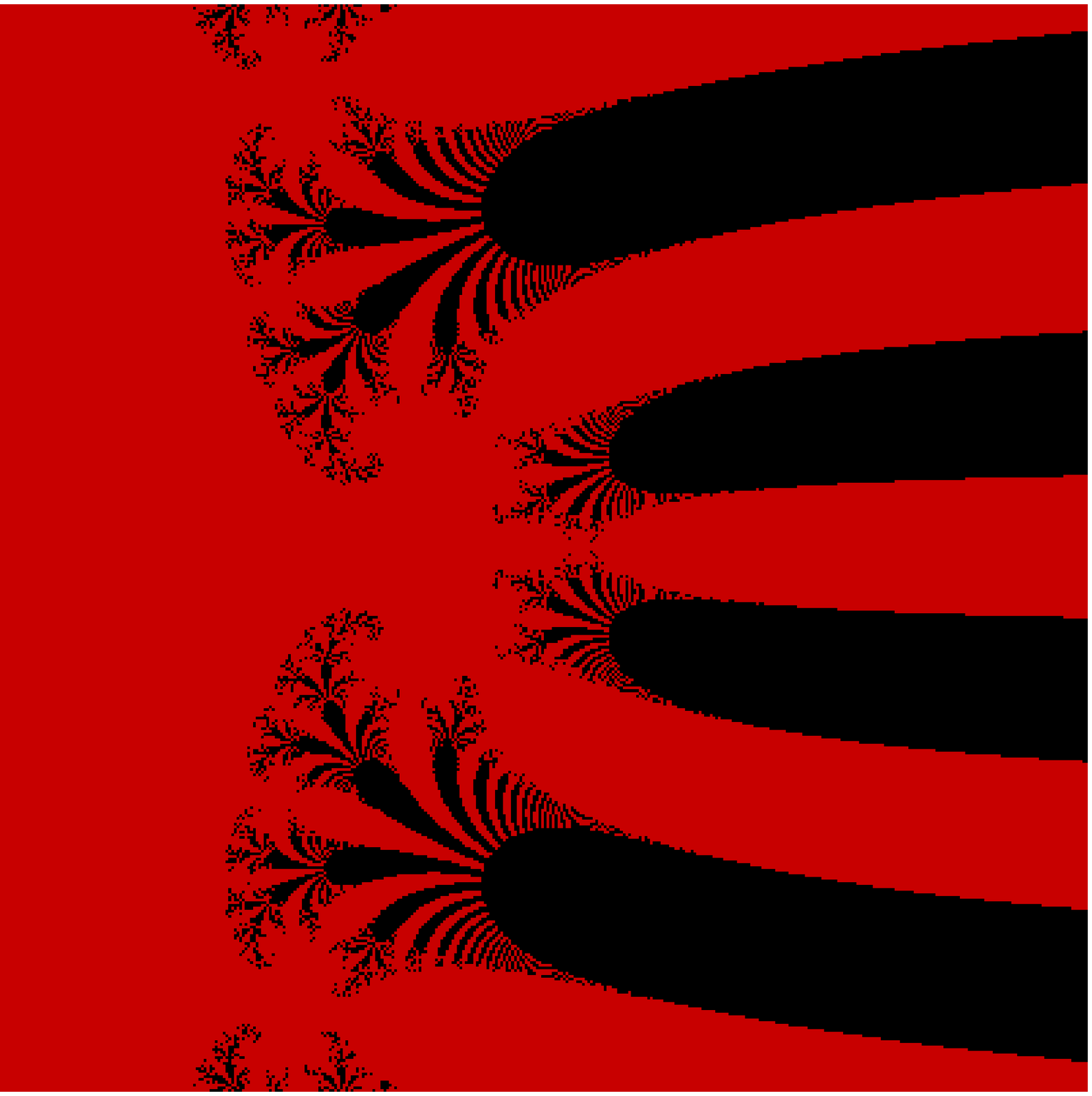}}
	\hspace{0.1in}
	\subfigure[\scriptsize{$\lambda = -8$}  ]{
	\includegraphics[width=0.3 \textwidth]{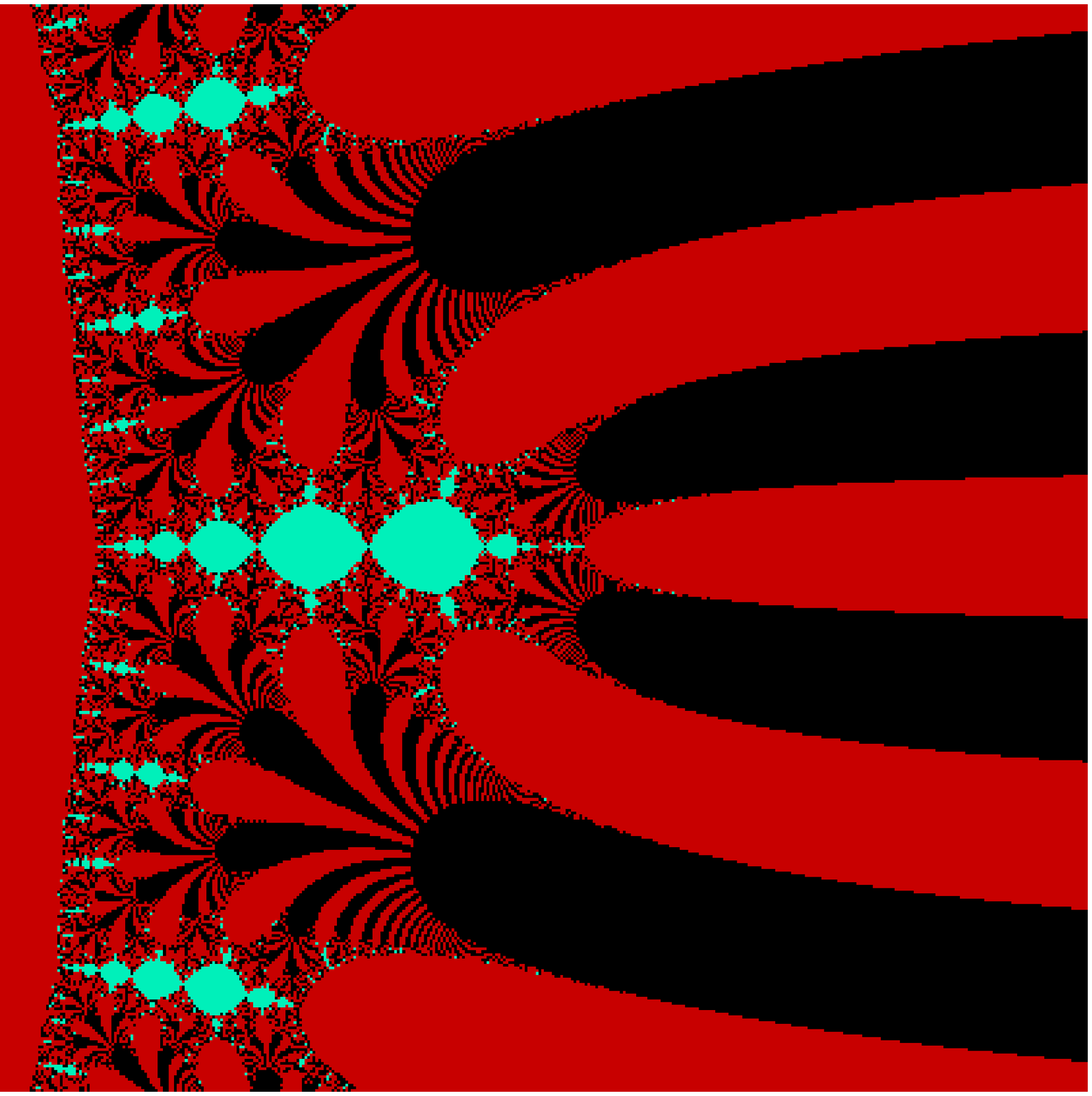}}
	\hspace{0.1in}
	\subfigure[\scriptsize{$\lambda = 6.9 $}  ]{
	 \includegraphics[width=0.3 \textwidth]{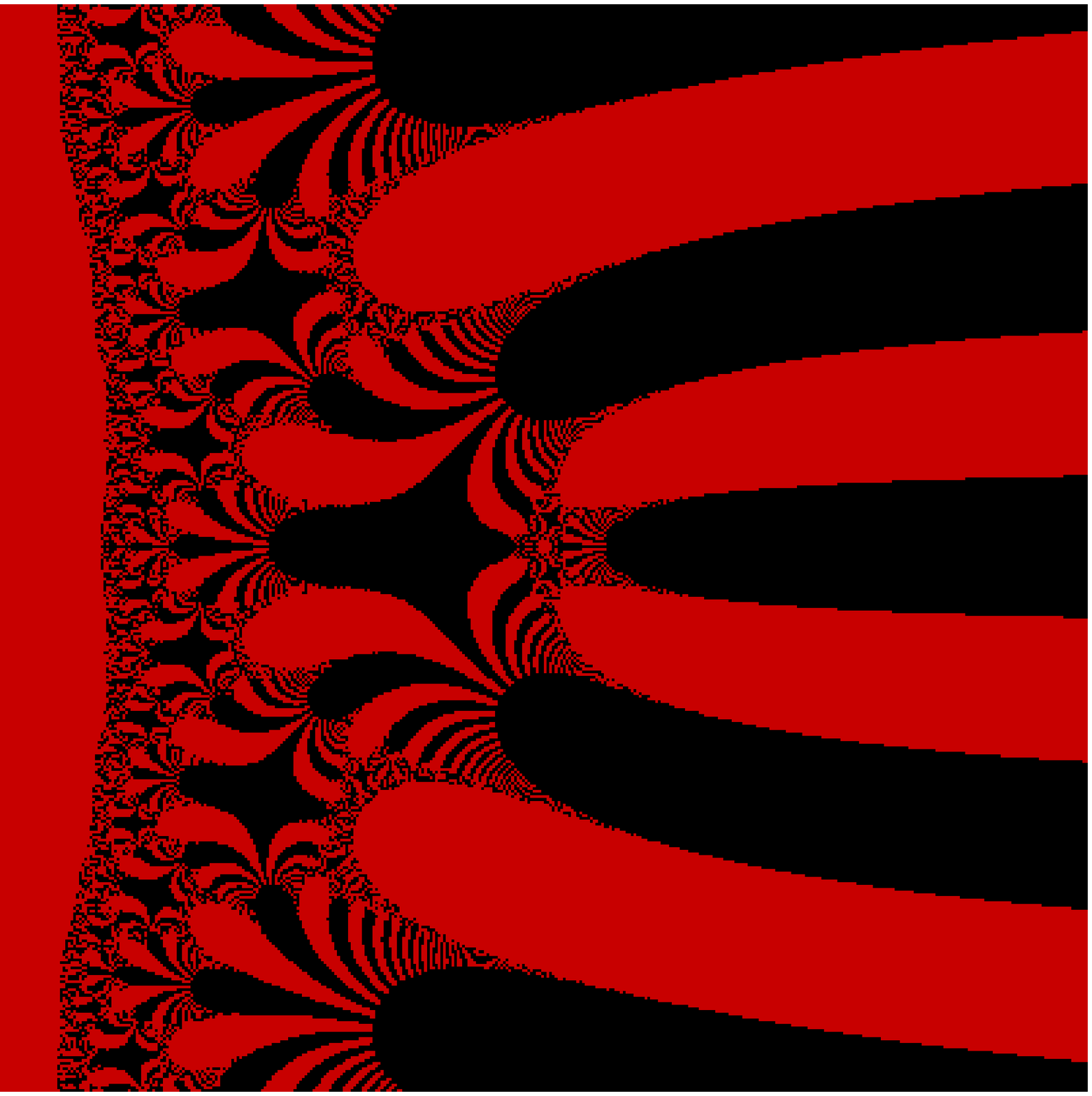}}
	\caption{\small{The Julia set for $F_{\lambda,2}$. Range $(-10,10) \times (-10,10)$.}}
	\label{jul2}
\end{figure}

\begin{figure}[hbt]
	\centering
	\subfigure[\scriptsize{ $m=2$} ]{
	 \includegraphics[width=0.3 \textwidth]{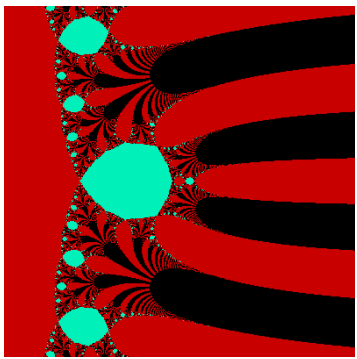}}
	\hspace{0.5in}
	\subfigure[\scriptsize{ $m=3$}  ]{
	\includegraphics[width=0.3 \textwidth]{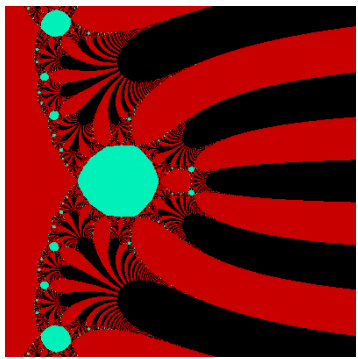}}
	\hspace{0.5in}
	\subfigure[\scriptsize{ $m=4$}  ]{
	 \includegraphics[width=0.3 \textwidth]{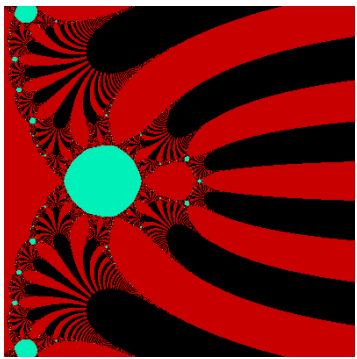}}
	\hspace{0.5in}
	\subfigure[\scriptsize{ $m=5$}  ]{
	 \includegraphics[width=0.3 \textwidth]{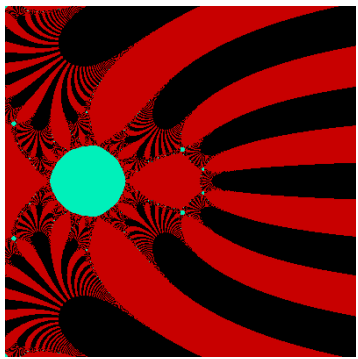}}
	\caption{\small{The dynamical plane of $F_{\lambda,m}$ for different values of $m$.
In every case  $\lambda=(-1)^{m-1} m (\frac{e}{m})^m$. Range $(-10,10) \times (-10,10)$.}}
	\label{julm}
\end{figure}

In Figs. \ref{jul2}-\ref{julm}, we display the Julia set of $F_{\lambda,m}$ for differents values of
$\lambda$ and $m$. The basin of attraction of $z=0$ is shown in red, while the components
of the Fatou set different from  A(0) are shown in blue. Points in the Julia set are shown in black.

In Fig. \ref{jul2}  we show the dynamical plane of function  $F_{\lambda,2} \, =
\, \lambda z^2 \exp(z) \,$, for three different values of $\lambda$. Apparently, the basin of $0$
contains an infinite number of horizontal strips, that extend to $+\infty$ as their real parts
tend to $+\infty$. Between these strips we find the well known structures, named Cantor Bouquets which
are invariant sets of curves governed by some symbolic dynamics. The existence of this kind of structures
in the Julia set are typical for critically finite entire transcendental functions (\cite{DT}).

As we change the parameter $\lambda$ we observe that the relative position of these bands
also changes, but not their width. Also, we can see the existence
of an unbounded region that extends to  $-\infty$ contained in $A(0)$.

In the next figure (Fig. \ref{julm}) we show a mosaic of different dynamical planes for some values of $m$,
specifically for $m=2, \, 3, \, 4,\,
\hbox{and } 5$. We choose $\lambda$ such that all these dynamical planes exhibit the same dynamical
behaviour, or more precisely, so that the critical point $z \, = \, -m$, is a superattracting
fixed point. A simple computation gives $\lambda \, = (-1)^{m-1} m (\frac{e}{m})^m$. In these
dynamical planes we see similar structures as in Fig. \ref{jul2}, even though the values of $m$ are
different.

We start with the following general result regarding $A(0)$.
\begin{lemma} \label{pa00}
$A(0)$ has either one or infinitely many connected components. Moreover, connected components different from $A^*(0)$ are
unbounded.
\end{lemma}
\begin{proof}
Using Sullivan's theorem (\cite{S}) we have that $A^*(0)$ is the unique fixed connected component of $A(0)$. For all other 
connected components of $A(0)$ there exists a number $ i>0$ such that $F_{\lambda,m}^i \,(U)\, = A^*(0)$, and $i$ is the smallest
number with this property. Suppose that there exist a finite number of connected components, and let $U_0 = A^*(0) \,
, \, U_1 \, , \, U_2 \, , \, ... \, , \, U_N$ the connected components of $A(0)$. We may choose the index $i$
in the natural way so that $F_{\lambda,m}^i \,(U_i)\, = A^*(0)$. Let $z \in U_N$ such that is not exceptional; then,
points in $F_{\lambda,m}^{-1}\, (z) \,$  belong  to $A(0)$, but not to $ U_0 \cup U_1 \cup ...
\cup U_N$, wich is a contradiction.

Now suppose that $U$ is a connected component of $A(0)$ different from $A^*(0)$, and let $i>0$ be the smallest number such
that $F_{\lambda,m}^i \,(U)\, = A^*(0)$. Let $z \in U$, and denote by $\gamma$ a simple path in $A^*(0)$ that
joins $F_{\lambda,m}^i(z)$ and 0. The preimage of $\gamma$ in $U$ must include a path $\gamma_1$ that joins $z$ and
$\infty$, since $0$ is an asymptotical value with no other finite preimage than itself. Thus we conclude that $U$ is unbounded.
\end{proof}

\subsection{Proof of Theorem A}

In this Sec. we describe the basin of attraction of the superattracting fixed point $z=0$.
Since  $z\, =\,0$ is a superattracting fixed point, there exists $\epsilon_0 \, > \, 0$ such that the
open disk  $D_{\epsilon_{0}} \, = \, \{ z \in {\mathbb C} \, ; \,|z| \, < \, \epsilon_0  \} $ is contained
in the immediate basin of attraction of $z=0$. First, we give an estimate of the size of the immediate basin of 
attraction, $A^*(0)$ (Proposition \ref{p21}), which will prove the first part of theorem A. Secondly, we find 
the first preimage of $D_{\epsilon_0}$
(Proposition \ref{hoho} and Proposition \ref{th}), proving the second part of theorem A. Finally, we find
the second preimage of $D_{\epsilon_0}$ (Proposition \ref{thbandas}), and prove the third part of
theorem A.

\begin{figure}[hbt]
\psfrag{m}[][]{\scriptsize $-m$}
\psfrag{fm}[][]{\scriptsize $h(-m)$}
\centerline{\includegraphics[width=0.8 \textwidth]{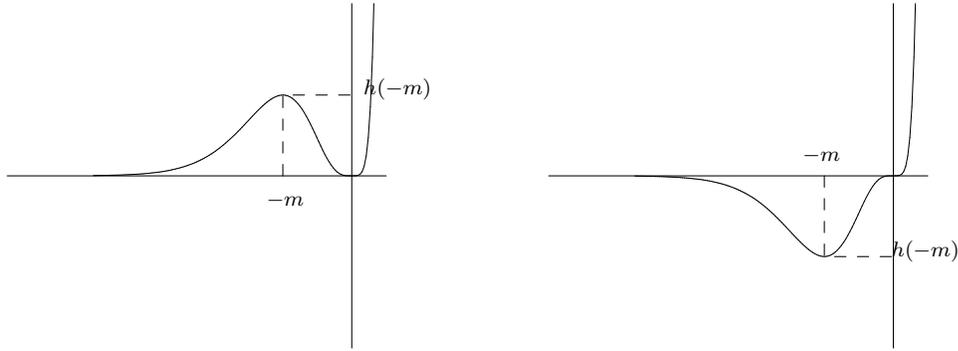}}
\caption{\small{Graph of $h(x)=x^m \exp(x)$. The left hand side corresponds to $m$ even and
 the right hand side to $m$
odd.}}
\label{hx}
\end{figure}

Before proving Proposition \ref{p21} we first look at some properties of the real funtion $h(x) \, = \, x^m e^x \, $
where $\, m \ge 1$. In Fig. \ref{hx} we show the graph of this function. It has a relative extremum
at $x=-m$, it is a monotone function on $(-\infty,-m)$ and it satisfies that
$|h(x)| \le |h(-m)| \, $ for all $\, x \le 0$. Also, it is an increasing function in $(0,+\infty)$.

Using these properties it is easy to prove the following auxiliary result.

\begin{lemma} \label{lem2}
Given  $r  \in (0,|h(-m)|]\,=\,(0,(\frac{m}{e})^m]$, the equation $ |h(x)| = r $ has a unique solution in $(-\infty,-m]$.
Moreover, given $s > 0$, the equation $h(x)=s$ has a unique solution in $(0,+\infty)$.
\end{lemma}

We now turn to estimate the function  $\epsilon_0(|\lambda|,m)$, and find its dependency on $\lambda \in  
{\mathbb C}-\{0\} \, \hbox{ and } m \ge 2$.

\begin{proposition}\label{p21}
If we define $\epsilon_0$ as the unique positive solution of the real
equation \par $x^{m-1}e^x = \frac{1}{|\lambda|}$, then $D_{\epsilon_0} \, = \, \{ z \in {\mathbb C} \, ; 
\, |z| \, < \, \epsilon_0 \} \subset A^*(0)$. Moreover, if $\lambda \in \br^{+}$ we have that $\epsilon_0$ lies 
in $\partial A^*(0)$.
\end{proposition}

\begin{proof}

In order to prove that  $D_{\epsilon_0}$ is contained in  $A^*(0)$, we use Schwartz's lemma. That is, it suffices to prove that
if $ |z| \, \le \, \epsilon_0  \,$ then $\, |F_{\lambda,m} (z)| \, \le \, \epsilon_0$.

Suppose $|z| \le  \epsilon_0$, we have

\[
|F_{\lambda,m} \, (z)  | \, =  \, | \lambda \,  z^m \, e^z |  \, = \,
| \lambda |  \, |z|^m e^{Re(z)} \, \le \, |z|\, |\lambda| \, |z|^{m-1} \, e^{|z|} .
\]

Since $|z|<\epsilon_0$ it follows that $|z|^{m-1} e^{|z|} < \epsilon_0^{m-1} e^{\epsilon_0}$, and using that
$\epsilon_0^{m-1} e^{\epsilon_0} = 1/|\lambda|$, we conclude
\[
|F_{\lambda,m} \, (z)  | \, \le \,
 |z| \, |\lambda| |z|^{m-1} e^{|z|} \, \le \, |z| \le \epsilon_0 .
\]

If $\lambda \in \br^{+}$, we have that $\lambda \epsilon_0^m \exp(\epsilon_0) \, = \epsilon_0$, i.e., $\epsilon_0$ is a 
fixed point. The multiplier of this fixed point is $\epsilon_0 + m >1$, and hence $\epsilon_0$ lies in the Julia set. By 
definition we have that $D_{\epsilon_0} \subset A^*(0)$, then $\epsilon_0$ lies in the boundary of $A^*(0)$.

\end{proof}

In the following auxiliary result we find a lower bound for $\epsilon_0$, which will be used in
the next Sec..

\begin{lemma}\label{lem_e0}
The value of $\epsilon_0 $ is always larger or equal than $ \min\{1,(\frac{1}{|\lambda|e})^{\frac{1}{m-1}}\}$
\end{lemma}

\begin{proof}

Suppose  $|\lambda| \, \ge \, 1/e$. This condition is equivalent to
$\frac{1}{|\lambda|e} \,\le \, 1$, hence we must prove that $\epsilon_0  \, \ge \,
(\frac{1}{|\lambda| \, e })^{1/(m-1)}$. Using that $x^{m-1}e^x$ is an increasing function on
$(0,+\infty)$, this condition is equivalent to
\[
\epsilon_0^{m-1} e^{\epsilon_0} \, \ge \, \frac{1}{|\lambda|e}\quad e^{(\frac{1}{|\lambda|e})^{\frac{1}{m-1}}}.
\]
By definition we have that $\epsilon_0^{m-1} e^{\epsilon_0} \, = \, \frac{1}{|\lambda|}$, then
\[
\frac{1}{|\lambda|} \, \ge  \,\frac{1}{|\lambda|e} \, e^{(\frac{1}{|\lambda|e})^{\frac{1}{m-1}}}
\]
or equivalently
\[
e \, \ge \,e^{(\frac{1}{|\lambda|e})^{\frac{1}{m-1}}}
\]

\noindent and this follows if $|\lambda| \, \ge  \,1/e $.

 If $|\lambda| \, \le \, 1/e$, we must prove that  $\epsilon_0 \, \ge  \,1$.
Using the same argument, i.e., that $x^{m-1} e^x$ is an increasing function, it follows that this condition
is equivalent to
\[
\epsilon_0^{m-1} e^{\epsilon_0} \, \ge e
\]
and this follows if $|\lambda| \, \le  \,1/e $.
\end{proof}

\begin{figure}[hbt]
\psfrag{HH}[][]{\small $H_{|\lambda|,m}$}
\psfrag{D}[][]{\small $D_{\epsilon_0}$}
\psfrag{x0}[][]{\small $x_0$}
\psfrag{yc}[][]{\small $y = C(x)$}
\psfrag{ymc}[][]{\small $y=-C(x)$}
\psfrag{x}[][]{\small $x$}
\centerline{\includegraphics[width=0.4 \textwidth]{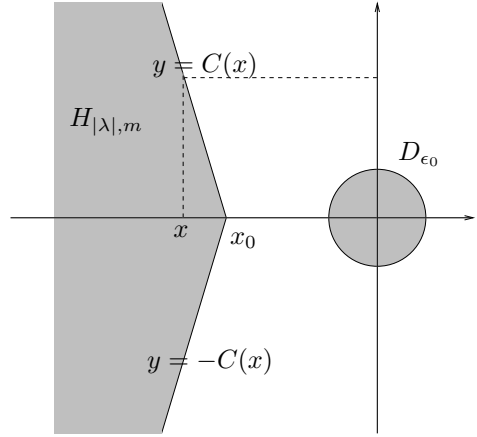}}
\caption{\small{Sketch of $H_{|\lambda|,m} \,.\, \hbox{ satisfying }  F_{\lambda,m} \, (H_{|\lambda|,m}) \,
\subset	\, D_{\epsilon_0} \subset A^*(0)$.}}
\label{hll}
\end{figure}

Next we want to find an open set $H_{|\lambda|,m} \, \subset \, {\mathbb C} $
such that $ F_{\lambda,m} (H_{|\lambda|,m} ) \, \subset  \, D_{\epsilon_0}\,$ (Fig. \ref{hll}). To that end, we first
obtain a value in ${\mathbb R}^{-}$, namely $x_0=x_0(|\lambda|,m)$, such that for all  $
x \in {\mathbb R}^{-}  \, \hbox{ with } \, x \le x_0 \,$ we have that  $\,  F_{\lambda,m}
(x) \in D_{\epsilon_0}$ (Proposition \ref{hoho}). After finding  $x_0$, we will look for an upper bound $C (x) \, \ge 0$,
such that if $z= x+yi, \, \hbox{with }\, x<x_0 \hbox{ and } |y|\le C(x) \, $, then $\,  F_{\lambda,m} (z) \, \in \,
D_{\epsilon_0}$ (Proposition \ref{th}).

\begin{proposition}\label{hoho}
 For all  $\lambda \in {\mathbb C} \hbox{ and } m \ge 2 \,$,
there exists  $x_0 \in (-\infty,-m]$, such that for all $ x \le x_0$, we have $\, F_{\lambda,m}(x) \in D_{\epsilon_0}$.
\end{proposition}
\begin{proof}
We suppose that $z=x+0i$ and we impose $|F_{\lambda,m}(z)| = \epsilon_0$, that is

\[
|F_{\lambda,m}(z)| \, = \, |\lambda| |x|^m e^x \, = \, \epsilon_0 \,
\]
\noindent or equivalently
\[
|h(x)|\,=\,|x|^m e^x \, = \, \frac{\epsilon_0}{|\lambda|},
\]
\noindent where $h(x)$ is the auxiliary function defined above.

If $|h(-m)| \le \frac{\epsilon_0}{|\lambda|}$, then we take  $x_0 = -m$, and for all $ x \in (-\infty,x_0) \,$,
we have $\, |h(x)| \le |h(-m)| \le \frac{\epsilon_0}{|\lambda|}$.
On the other hand, if $|h(-m)| > \frac{\epsilon_0}{|\lambda|}>0$,
we define $x_0$ as the unique solution of equation $|h(x)|=\frac{\epsilon_0}{|\lambda|}$ in the interval $(-\infty,-m)$.
Since  $|h(x)|$ is an increasing function in $(-\infty,-m)$, it follows that for all
$ x \in (-\infty,x_0) \,$, then $\, |h(x)| \le |h(x_0)| =  \frac{\epsilon_0}
{|\lambda|}$.
\end{proof}

\begin{proposition} \label{th}
Let $x_0=x_0(|\la|,m)\,$ be as in Proposition \ref{hoho}. There exists $\, C(x) \ge 0$ such that the 
open set (Fig. \ref{hll})
 $$ H_{ |\lambda|,m} \, = \, \left\{ z \, =
\, x \, + \, y i \, \left| \begin{array}{ll}
x \,&  \in \, (-\infty,x_0)  \\
y \,&  \in  (-C(x),C(x))  \end{array}\right. \right\}  $$

satisfies  $F_{\lambda,m} \, (H_{|\lambda|,m}) \, \subset \, D_{\epsilon_0}$.
\end{proposition}

\begin{proof}
We want to calculate an upper bound  $C(x)  \ge  0$, such that if $z \, = \, x +
yi \, $, with $x \in (-\infty,x_0) \, \hbox{ and } \, |y| \le C(x) \, $, then $\,
F_{\lambda,m}(z) \in D_{\epsilon_0}$. If we require that  $F_{\lambda,m}(z) \, \in \,  D_{\epsilon_0}$,
we obtain the definition of $C(x)$. Let $z \, = \, x \pm y i $, with  $x \in (-\infty,x_0)$. Then

\[
|F_{\lambda,m} (z)| \, = \,  |\lambda| |z|^m \exp(Re(z)) \, = \,|\lambda| \left[ \sqrt{x^2  + y^2} \right]^m \exp(x) \, = \, |\lambda| (x^2 \, + \, y^2)^{m/2} \exp(x).
\]

\noindent Using the expression above, and requiring  $|F_{\lambda,m}(z)| \, \le \,  \epsilon_0 $, we obtain

\[
|y| \, \le \,
+ \sqrt{ \,\left[ \frac{\epsilon_0}{|\lambda|} \right]^{2/m} \exp(-x\frac{2}{m})\, -
\, x^2 }.
\]

\noindent  Thus, the right hand side of this inequality gives an analytic expression for the function $C(x)$.
\end{proof}

\begin{figure}[hbt]
\psfrag{a}[][]{\scriptsize $y=A|x|$}
\psfrag{b}[][]{\scriptsize $y=A|x|^2$}
\psfrag{c}[][]{\scriptsize $H_{|\lambda|,m}$}
\centerline{\includegraphics[width=0.4 \textwidth]{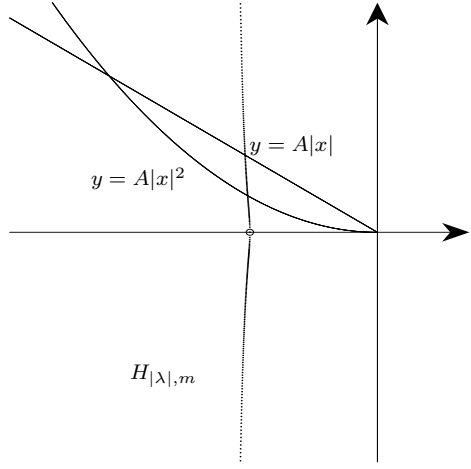}}
\caption{\small{Relation between $ H_{|\lambda|,m}\, \hbox{ and } \, y=A|x|^k \, \hbox{ for } \, k=1,2$}}
\label{esq}
\end{figure}

The next lemma gives a simpler condition to assure that a point
 $z \in {\mathbb C}$ lies in $H_{|\lambda|,m}$.
See Fig. \ref{esq}.
\begin{lemma} \label{pro}
The point $z= x +  y i \, $ lies in $ H_{|\lambda|,m}$ if there exists  $ \, k \ge 1 \, \hbox{ and } \, A \ge 0$, such 
that  $|y| \, \le \, A\,|x|^k$ and $|x| \,$ is large enough.

\end{lemma}

\begin{proof}
We will prove that
$A|x|^k \, \le \, C(x) \hbox{ as } x \to -\infty$. Using the definition of $C(x)$, this is equivalent to showing

\[
A|x|^k \, < \, \sqrt{ \,\left[ \frac{\epsilon_0}{|\lambda|} \right]^{m/2}
\exp(-x\frac{m}{2})\, - \, x^2 },
\]
\noindent or

\[
\left[A^2|x|^{2k} \, + \, x^2 \right]\exp(x\frac{m}{2}) \, < \,
\left[ \frac{\epsilon_0}{|\lambda|}\right]^{m/2}.
\]

\noindent The left hand side of this inequality is a function that tends to zero as $x \, $ tends to $\, -\infty$, whereas
the right hand side is positive.

\end{proof}

We proceed now to the third iterate, by proving the existence of some strips in dynamical plane, such that
the image of these open sets under $F_{\la,m}$ is contained in $H_{|\lambda|,m}$ (see Fig. \ref{a0}).

Before calculating the preimage of the set $H_{|\lambda|,m}$, we first find the preimage of
the negative real axis under the function $F_{\lambda,m}$. Hereafter, we denote by $Arg(.) \in (-\pi,\pi]$ the
argument. Using the definition of  $F_{\lambda,m}$ it is easy to see that

\[
Arg(F_{\lambda,m}(z)) = Arg(\lambda) + m Arg(z) + Im(z) \qquad (\hbox{ mod } 2 \pi).
\]

Finding the preimages of  ${\mathbb R}^-$ is equivalent to solving
\[
Arg(F_{\lambda,m}(z)) = \pi.
\]
We denote $r=|z|$ and $\alpha = Arg(z)$. Then the equation above is equivalent to
\[
Arg(\lambda) + m \alpha + r sin(\alpha) = (2k+1)\pi \qquad k \in {\mathbb Z}.
\]
Hence, we obtain
\[
r = \rho(\alpha) =  \frac{(2k+1)\pi - m \alpha - Arg(\lambda)}{sin(\alpha)} \qquad \alpha \in (-\pi,\pi).
\]
We denote each of these curves by $\sigma_k\, =\,\sigma_k(\lambda,m)$, where the possible values of the argument depend on $k$. More precisely,

if  $m=2j$ for $j \in \bz$
\[
\sigma_k = \rho(\alpha) e^{i\alpha} \hbox{ where }    \left\{  \begin{array}{llll}    0< \alpha < \pi \qquad & \hbox{  if  }\quad  k \ge j \\
  0< \alpha <\frac{(2k+1)\pi-Arg(\lambda)}{m} \qquad & \hbox{  if  }\quad 0 \le  k \le j-1 \\    \frac{(2k+1)\pi-Arg(\lambda)}{m}< \alpha <0 \qquad & \hbox{  if  }\quad -j \le  k \le 0 \\     -\pi< \alpha <0  \qquad & \hbox{  if  }\quad k \le -(j+1) \\ \end{array}\right. ;
\]

if  $m=2j+1$ for $j \in \bz$
\[
\sigma_k =  \rho(\alpha) e^{i\alpha} \hbox{ where } \left\{  \begin{array}{llll}   0< \alpha < \pi \qquad & \hbox{  if  }\quad  k \ge j+1 \\
0< \alpha <\frac{(2k+1)\pi-Arg(\lambda)}{m} \qquad & \hbox{  if  }\quad  0 \le  k \le j \\
 \frac{(2k+1)\pi-Arg(\lambda)}{m}< \alpha <0 \qquad & \hbox{  if  }\quad -j \le  k \le 0 \\
-\pi< \alpha <0  \qquad & \hbox{  if  }\quad k \le -(j+1) \\ \end{array}\right.
\]

\begin{figure}[hbt]
	\centering
	\subfigure[\scriptsize{Graph of $\sigma_k$ for $\la = 0.45 + 0.35i$ and $m=5$.} ]{
\psfrag{a}[][]{\scriptsize $\sigma_0$}
\psfrag{b}[][]{\scriptsize $\sigma_1$}
\psfrag{c}[][]{\scriptsize $\sigma_2$}
\psfrag{d}[][]{\scriptsize $\sigma_{-1}$}
\psfrag{e}[][]{\scriptsize $\sigma_{-2}$}
\psfrag{m}[][]{\scriptsize $\sigma_{-3}$}
\psfrag{f}[][]{\scriptsize $\sigma_{-4}$}
\psfrag{g}[][]{\scriptsize $\sigma_{-5}$}
\psfrag{h}[][]{\scriptsize $\sigma_{-6}$}
\psfrag{i}[][]{\scriptsize $\sigma_3$}
\psfrag{j}[][]{\scriptsize $\sigma_4$}
\psfrag{k}[][]{\scriptsize $\sigma_5$}
\psfrag{l}[][]{\scriptsize $\sigma_6$}
 \fbox{\includegraphics[width=0.4 \textwidth]{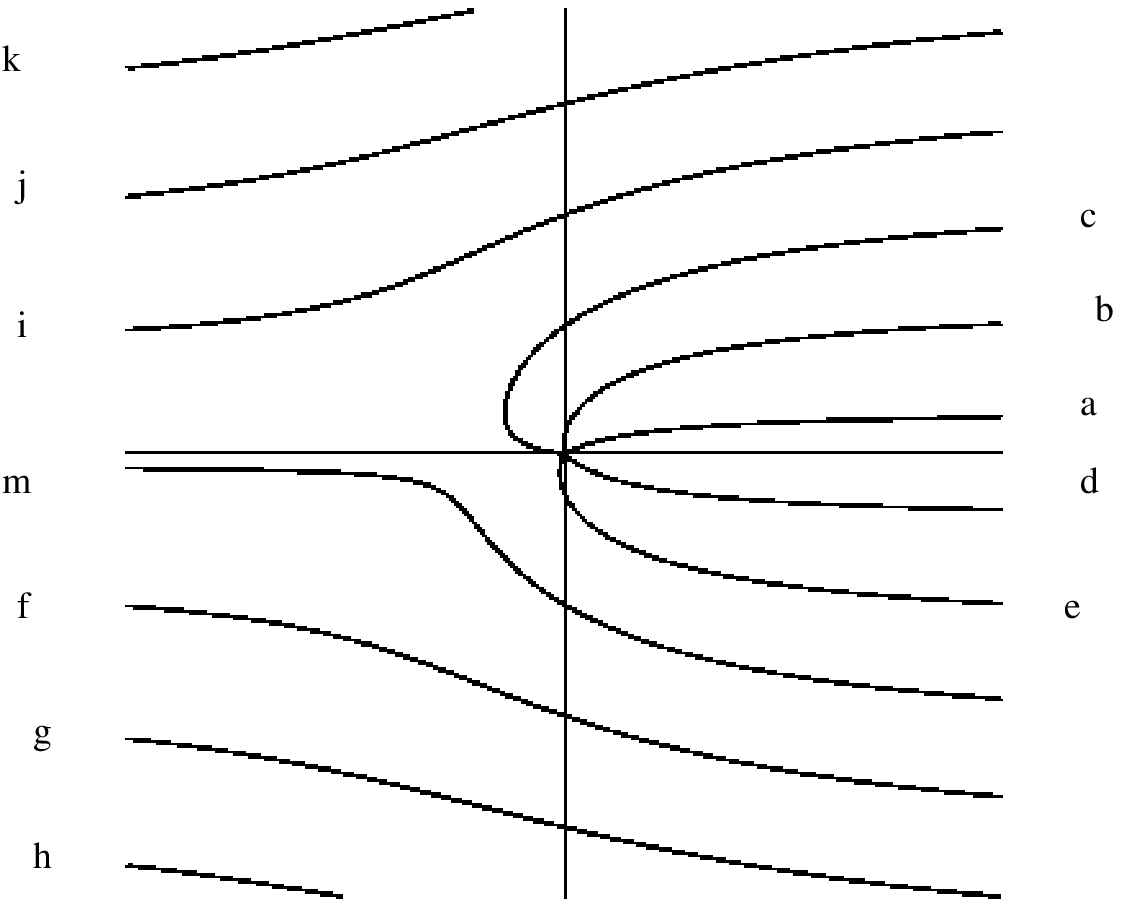}}}
	\hspace{0.5in}
\subfigure[\scriptsize{The julia set of $F_{0.45+0.35i,5}$} ]{
\includegraphics[width=0.4 \textwidth]{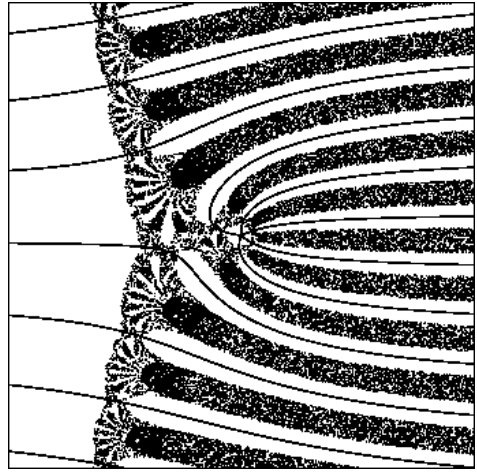}}
	\caption{\small{Strips in the dynamical plane.}}
	\label{par2}
\end{figure}

In Fig. \ref{par2} we show some of these curves for $m=5$. As their real parts tend to $+\infty$, the
$\sigma_k$'s are asymptotic to the lines
$Im(z) = (2k+1)\pi-Arg(\lambda)$. There are $m$ of these curves, that start at the origin
and tend to $+\infty$. The others are asymptotic to the lines
$Im(z) = (2k+1)\pi-Arg(\lambda)\,$  when $k<0$, or  $Im(z) = 2k\pi-Arg(\lambda)\,$ when $k>0$.

Now we may find preimages of the open set  $H_{|\lambda|,m}$. First we find preimages of the interval
$(-\infty,x_0)$. There exists one preimage of this interval on each curve $\sigma_k$. Moreover, the preimage of a real
number tending to $-\infty$, is a complex number on $\sigma_k$ whose real part tends to $+\infty$. Hence, the preimages of the set $H_{|\lambda|,m}$ contain some strips, namely  $S_{\lambda,m}^k$, around $\sigma_k$. (See Fig. \ref{par2}).

Let  $z \in \sigma_k$. If we evaluate $F_{\la,m}(z)$, we have
\[
F_{\la,m}(z)= |\la| |z|^m e^{Re(z)} \, e^{(Arg(\la)+m Arg(z) + Im(z))i} = - |\la| |z|^m e^{Re(z)}
\]

\noindent since each  $\sigma_k$ is a preimage of the negative real axis. The expression above shows that
if we keep $Re(z)$  constant, and we increase the index $k$ of the curve $\sigma_k$, we obtain values tending to
 $-\infty$, since $|z|$ increases. Hence, if we denote by $q_0(k)$ the preimage of $x_0$ on $\sigma_k$,
its real part decreases as $|k|$ increases (at least after a certain point). This fact explains
the apparent arangement of the strips in dynamical plane (Fig. \ref{par2}).

We now  prove that these strips have an asymptotic
width equal to $\pi$. We fix a value $k \in {\mathbb Z}$, and we recall that $\sigma_k$ tends asymptotically to the line
 $Im(z) = (2k+1) \pi - Arg( \lambda)$, as its real part tends to $+\infty$.

\begin{proposition} \label{thbandas}
Given any $\, y \in ( (2k+1) \pi - Arg( \lambda) - \frac{\pi}{2} \,  , \, (2k+1) \pi - Arg( \lambda) + \frac{\pi}{2}) , \, $ there exists
a real number $\, x_*$ such that for all $\,x \ge x_* $, $F_{\lambda,m} (x+yi) \, \in  \, H_{|\lambda|,m}$.

\end{proposition}

\begin{proof}

Let  $z = x + i y $, where $ y \in ( (2k+1) \pi - Arg( \lambda) - \frac{\pi}{2} \,  , \, (2k+1) \pi - Arg( \lambda) +
\frac{\pi}{2}) $. If we write the parameter $\la$ in polar
coordinates, $\lambda \, = \, se^{i\beta} $, then

\[
\begin{array}{ll}
F_{\lambda,m} (z) \, = \, & \lambda \, z^m \, \exp(z) \, = \, s e^{i\beta} \, (x \, + y i)^m
e^x \, e^{i y} \, = \, s (x \, + y i)^m \, e^x \, e^{i(y \, + \, \beta)} \, = \, \\
 & \{P(x) \, + Q(x)i \} \, e^x \, \{\cos(y \, + \, \beta) \, + \, \sin(y \, + \,
\beta)i \}
\end{array}
\]
\noindent where
\[
 \left\{  \begin{array}{ll} P(x) \, = & Re(s(x \, + \,y i)^m) \, = \,
sx^m \, + \, {\cal O} (x^{m-2}) \\
 Q(x) \, = & Im(s(x \, + \,y i)^m) \, = \,
msy x^{m-1} \, + \, {\cal O} (x^{m-3}) . \end{array}\right.
\]

\noindent Using the expressions above in $F_{\lambda,m} (z)\,$, we obtain

\begin{equation}
\begin{array}{llllll}
F_{\lambda,m} (z) &  =  \{P(x)\cos(y  +  \beta) - Q(x)\sin(y  +  \beta) \} e^x + \\
                     & \quad \{ P(x)\sin(y +  \beta) + Q(x)\cos(y  + \beta) \}e^x  i  \\
 & = \{ P(x)\cos(y \, + \, \beta) \, + {\cal O}(x^{m-1}) \} e^x \, + \, \\ \label{eqnum}
& \quad  \{ P(x)\sin(y \, + \, \beta) \, + {\cal O}(x^{m-1}) \} \, e^x \, i  \, \\
 & = \{sx^m\cos(y\, + \, \beta) \, + \, {\cal O}(x^{m-1}) \} e^x \,  + \, \\
& \quad  \{sx^m\sin(y \, + \, \beta) \, + \, {\cal O}(x^{m-1}) \} e^x \, i.
\end{array}
\end{equation}

\noindent We recall that $ \, y  \in  (\frac{\pi}{2}  +  2k\pi -  \beta \, ,
 \, \frac{3\pi}{2}  +  2k\pi  -  \beta ) \, $, or equivalently $ \, y \, + \, \beta \in (\frac{\pi}{2} \, + \, 2k\pi , \frac{3\pi}{2} \, + \, 2k\pi) \,$  which implies that $ -1 \le cos(y \, + \, \beta ) < 0$

Since $Re(F_{\lambda,m} (z)) \, = \, \{sx^m cos(y \, + \, \beta )  \, + \, {\cal O}
(x^{m-1}) \}e^x$, it follows that
\[
\lim_{x \, \to \, +\infty} Re(F_{\lambda,m} (z)) \, = \,
\lim_{x \, \to \, +\infty} \{ cos(y \, + \, \beta ) x^m s \, + \, {\cal O} (x^{m-1}) \}e^x \, = \,
 -\infty.
\]

Hence, there exists $x_*$ large enough, such that  $Re(F_{\lambda,m} (z)) \, < \, x_{0} $ for $x \ge x_*$.

 Finally, we use lemma \ref{pro} to conclude the proof. From equation (\ref{eqnum}) we have

\[
\lim_{x \to +\infty} \frac{Re(F_{\lambda,m}(z))}{Im(F_{\lambda,m}(z))} \, =  \,
\lim_{x \to +\infty} \frac{sx^m \cos(y \,+\beta) \, + \, {\cal O} (x^{m-1})}
{sx^m \sin(y \,+\beta) \, + \, {\cal O} (x^{m-1})} \, = \, \frac{\cos(y \,+\beta)}
{\sin(y \,+\beta)}.
\]

Hence, if $\sin(y \, + \, \beta) \neq 0 \,$, it follows that $\, Im(F_{\lambda,m}(z)) \, \approx \,
\frac {\sin(y \,+\beta)}{\cos(y  + \beta)} Re(F_{\lambda,m}(z)) $, and we can apply
lemma \ref{pro} ( with $A \, = \, |\frac{\sin(\alpha \,+\beta)}{\cos(\alpha  + \beta)}|
\, \hbox{ and } \, k=1$). Using this lemma, we conclude that  $ F_{\lambda,m}(z)\, $ lies
in  $\, H_{|\lambda|,m}$, if $x$ is large enough.

On the other hand, if  $\sin(y \, + \, \beta) = 0 \,,$ then  $\, cos(y \, + \, \beta) = -1$, and we obtain

\[
\begin{array}{ll}
F_{\lambda,m} (z) \,  & = \, \{-P(x) -Q(x)i\}e^x \, \\
  & = \{-sx^m +  {\cal O} (x^{m-1})\}e^x \, + \,
\{m s y x^{m-1}+  {\cal O} (x^{m-2})\}ie^x.
\end{array}
\]

\noindent Hence
\[
\lim_{x \to +\infty} \frac{Re(F_{\lambda,m}(z))}{Im(F_{\lambda,m}(z))} \, =  \,
 \lim_{x \to +\infty} \frac{-sx^m + {\cal O} (x^{m-1})}{-msyx^{m-1}+  {\cal O} (x^{m-2})} \, = \,
+\infty .
\]

If $x$ is large enough, there exists $K>0$ such that

\[
\left| \frac{Re(F_{\lambda,m}(z))}{Im(F_{\lambda,m}(z))} \right|
\, > K \,
\]

\noindent that is, $ | Im(F_{\lambda,m}(z)) | \,  \le \, \frac{1}{K} |Re(F_{\lambda,m}(z))|$. Hence, using lemma \ref{pro}
 ($A = \frac{1}{K}$ , $k=1$), we also obtain that $ F_{\lambda,m}(z)\,$ lies in $\,H_{|\lambda|,m}.$

\end{proof}


\section{The Parameter Planes}

The orbit of the free critical point $z=-m$, determines in large measure the dynamics of $F_{\la,m}$.
Indeed, the functions $F_{\la,m}(z) = \la z^m \exp(z) $ are entire maps with a finite
number of critical and asymptotic values. These kind of functions do not have wandering domains
nor Baker domains. By the Sullivan classification, we know that
if the orbit of $z=-m$ tends to $\infty$ then the Fatou set must coincide with the basin of $0$, i.e.,
${\cal F} (F_{\lambda,m}) \, = \, A(0)$, since no other Fatou components can exist besides those that belong to
$A(0)$. The set $B_m$ is defined as before as

$$B_m \, = \, \{ \lambda \in  {\mathbb C} \, | \, F_{\lambda,m}^{\circ n} (-m)  \nrightarrow
\infty \}.$$

\begin{figure}[hbt]
	\centering
	\subfigure[\scriptsize{Range $(-25,25) \times (-25,25)$} ]{
	\includegraphics[width=0.45 \textwidth]{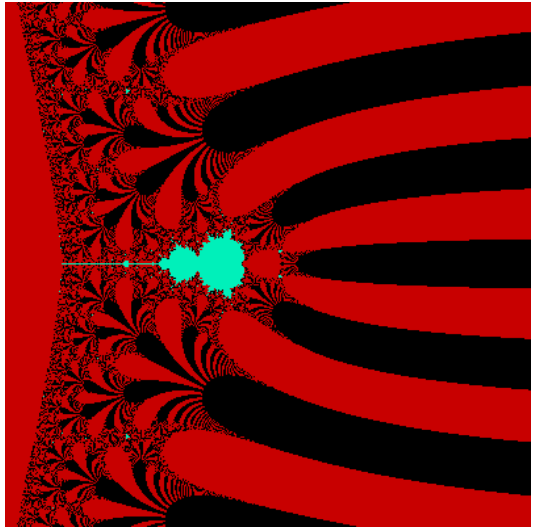}}
	\hspace{0.5in}
	\subfigure[\scriptsize{Range $(-15,5) \times (-8,8)$}  ]{
	\includegraphics[width=0.45 \textwidth]{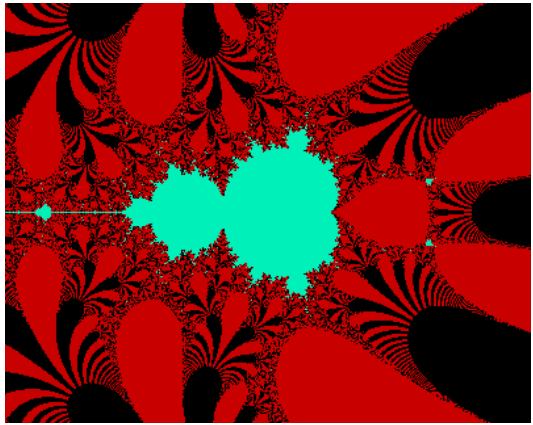}}
	\caption{\small{Parameter plane for $F_{\lambda,2}$. Color codes are explained in the text.}}
	\label{pp2}
\end{figure}

In each of these sets, we may also distinguish between two different behaviours: those parameters values for which
$-m \in A(0)$ and those for which this does not occur. Let $\stackrel{\circ}{B_m}$ denote the interior of $B_m$.

\begin{definition} Let $ U $ be a connected component of $\stackrel{\circ}{B_m}$. We say that $U$ is a {\it capture zone} if 
for all $\, \lambda $ in  $\, U$ it is true that $\,\lim_{n \to +\infty} F_{\lambda,m}^{\circ n}(-m) = 0$, or in other words,
 $ -m \in A(0)$. We then say that the orbit of the critical
point is captured by the basin of attraction of the superattracting fixed point $z=0$.
\end{definition}

\begin{figure}[hbt]
	\centering
	\subfigure[\scriptsize{$B_3$. Range $(-12,12) \times (-12,12)$} ]{
	 \includegraphics[width=0.36 \textwidth]{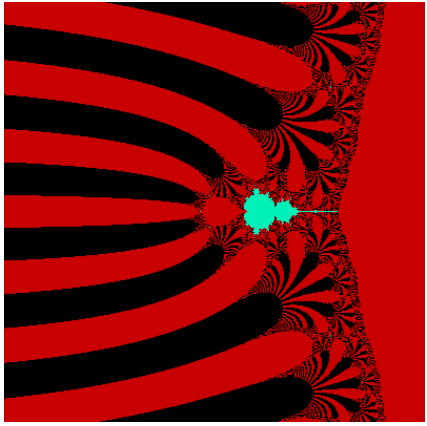}}
	\hspace{0.5in}
	\subfigure[\scriptsize{$B_4$. Range $(-4,2) \times (-3,3)$}  ]{
	 \includegraphics[width=0.36 \textwidth]{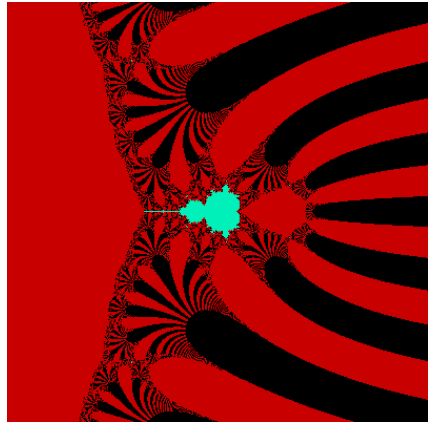}}
	\hspace{0.5in}
	\subfigure[\scriptsize{$B_5$. Range $(-0.8,0.8) \times (-0.8,0.8)$}  ]{
	 \includegraphics[width=0.36 \textwidth]{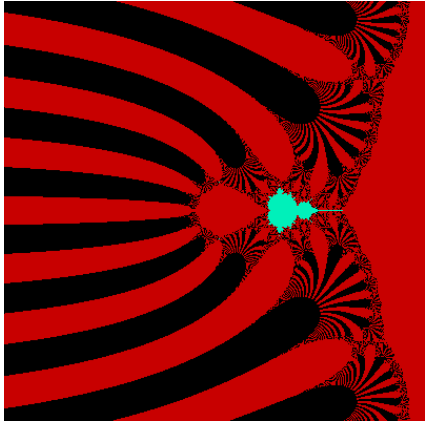}}
	\hspace{0.5in}
	\subfigure[\scriptsize{$B_6$. Range $(-0.15,0.15) \times (-0.15,0.15)$}  ]{
	 \includegraphics[width=0.36 \textwidth]{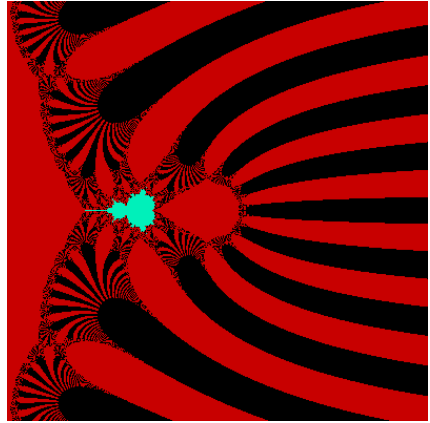}}
	\caption{\small{Parameter plane for $F_{\lambda,m}$, for differents values of $m$.}}
	\label{ppm}
\end{figure}

In Figs. \ref{pp2}- \ref{ppm}, we show a numerical approximation of the set $B_m$ for different values of $m$.
The capture zones are shown in red, while other components of $B_m$  are
shown in blue. The parameter values for which the orbit of the free critical point tends to $\infty$ are shown
in black. In these sets we can see a countable quantity of horizontal strips. If $m$ is even these strips
extend to $+\infty$ as
the real part of $\la$ tends to $+\infty$, whereas if $m$ is odd these strips
extend to $-\infty$ as the real part of $\la$ tends to $-\infty$. Notsurpringly, the distribution of these capture zones
in the parameter plane (Fig. \ref{zona}) appears to be similar to the distribution of $A(0)$ in the dynamical plane (Fig. \ref{a0}).

We start with the following simple facts.
\begin{proposition}
Let $U$ be a capture zone of $B_m$ and let $ \lambda \in U\,$. Then $ \,{\cal F} (F_{\lambda,m}) \, =
\, A(0)\,$ and $\,{\cal J} (F_{\lambda,m}) \, = \, \partial A(0)$.
\end{proposition}

\begin{proof}
As in the case of the critical value tending to $\infty$, since the only free critical point of $F_{\la,m}$ lies
in the basin of $0$, no other components different from those in $A(0)$ can exist in ${\cal F} (F_{\lambda,m})$.
Let ${\cal D}$ be the union of all the components of the Fatou set, then ${\cal J} \, (F_{\lambda,m}) \, = \, \partial
{\cal D}$ (\cite{CG}). In our case, if  $\lambda \in U \, $ and  $\, U$ is a capture zone, then  ${\cal D} \,
 = \, A(0).$
\end{proof}

The main objective of this Sec. is to describe the most obvious capture zones contained in
$B_m$, as well as to describe the dynamical plane for parameter values that belong to such components.

\subsection{Proof of Theorem B}

In this Sec. we describe the main capture zone $C_m^0$. We recall that
\[
C_m^n = \{ \lambda \in \stackrel{\circ}{B_m} | F_{\lambda,m}^n (-m) \in A^*(0) \hbox{ and $n$ is the smallest number with this property} \}
\]
We prove each statement of theorem B in a different proposition.

\begin{proposition}\label{ant}
The critical point $-m$ belongs to $A^*(0)$ if and only if the critical value $F_{\lambda,m}(-m)$ belongs 
to $A^*(0)$. Hence $C_m^1 \, = \, \o$.
\end{proposition}
\begin{proof}
Supose that $F_{\lambda,m}(-m) \in A^*(0)$. Let $\gamma$ be a simple path in $A^*(0)$ that joins $F_{\lambda,m}(-m)$ and $0$. The
set of preimages of $\gamma$  must include a path $\gamma_1$ that joins $-\infty$ with $-m$, and also a path $\gamma_2$ that
joins $-m$ and $0$ (since $-m$ is a critical point and $0$ is a fixed point and asymptotic value). Hence $\gamma_1 \cup 
\gamma_2 \subset A^*(0)$ and so does $-m$. Conversely, if $-m \in A^*(0)$ we have that $F_{\lambda,m}(-m) \in A^*(0)$.
\end{proof}

\begin{proposition} \label{pepe}
The set $C_m^0$ contains the disk $\{ \lambda \in {\mathbb C} \, ;\, |\lambda| \, < \,
\min(\frac{1}{e},(\frac{e}{m})^m ) \}$.
\end{proposition}
\begin{proof}
We denote  $D_m$ the open disk  $ \{ \lambda \in {\mathbb C} \, ;\, |\lambda| \, < \,
\min(\frac{1}{e},(\frac{e}{m})^m ) \} \,$. Let $\lambda \in D_m$, we will prove that $F_{\lambda,m}(-m) $
lies in $D_{\epsilon_0}$ which we know belongs to $A^*(0)$. In order to do so, we use that
$\epsilon_0 \ge \min(1,(\frac{1}{|\lambda| e})^{1/(m-1)})$ (lemma \ref{lem_e0}).
We choose $\lambda \in D_m$. Then  $\, |\lambda| < \frac{1}{e} \,$, and hence
$\, \epsilon_0 \ge 1$. The condition $\la \in D_m$ also implies that
 $|\lambda| <  (\frac{e}{m})^m$. Hence
\[
|F_{\lambda,m}(-m)| \, = \, |\lambda| |(-m)^m e^{-m}| \,= \,
 |\lambda| \left(\frac{m}{e}\right)^m \, < \, 1
\le \epsilon_0,
\]
\noindent and $F_{\lambda,m}(-m) \, $ lies in $ \,A^*(0)$.
\end{proof}

\begin{proposition}\label{prim}
The set $C_m^0$ is bounded. In fact it is contained in the closed
disk $\{ \lambda \in {\mathbb C} \, ;\, |\lambda|
\le (\frac{e}{m-1})^{m-1} \}$.
\end{proposition}

\begin{proof}
We will prove that $-m \notin A^*(0)$ for all $\lambda \in {\mathbb C} $ such that $|\lambda| > (\frac{e}{m-1})^{m-1}$.
Let  $D$ the disk centered at $0$ of radius $m-1$. If we calculate the image of its
boundary, $\{|z|=m-1\}$, we obtain

\[
|F_{\lambda,m}(z)|= |\lambda| |z|^m e^{Re(z)} \ge  |\lambda| (m-1)^m e^{-(m-1)} > m-1
\]

\noindent where the inequality is obtained using $|\lambda| > (\frac{e}{m-1})^{m-1}$. This shows that $D \subset F_{\lambda,m}(D)$, and
hence $A^*(0) \subset D$. Since $-m \notin D$, the proposition follows.
\end{proof}

\begin{proposition}
If $\lambda \in C_m^0$ then $A(0)=A^*(0)$, i.e., the basin of attraction of $z=0$ has a unique connected
component and hence it is totally invariant. Moreover, the boundary of $A^*(0)$ (which equals the Julia set) is a Cantor 
bouquet and hence it is disconnected and non-locally connected.
\end{proposition}
\begin{proof}
Let $\lambda \in C_m^0$. As in proposition \ref{ant}, let $\gamma$ be a simple path in $A^*(0)$ that 
joins $F_{\lambda,m}(-m)$ and $0$. The
preimage of $\gamma$ must include a path $\tilde{\gamma}$ contained in $A^*(0)$ that joins $-\infty$ with $0$ passing through 
$-m$ ($\tilde{\gamma}$ maps 2-1 to $\gamma$). Since $H_{|\lambda|,m}$ intersects  $\tilde{\gamma}$ so it follows that $H_{|\lambda|,m} \subset A^*(0)$. All preimages of $\tilde{\gamma}$,
are contained in $A^*(0)$ as well, since they all intersect $H_{|\lambda|,m}$. In fact, we have
that $A(0)=A^*(0)$ since any preimage of $D_{\epsilon_0}$ must contain points of $H_{|\lambda|,m}$. Hence $A(0)$ has a unique
connected component. In fact, from \cite{DG},\cite{BD}, it follows that the Julia set has an uncountable number of connected
components
and it is not locally connected at any point.

Using \cite{DT}  one can show that the Julia set contains a
Cantor Bouquet tending to $\infty$ in the direction of the positive real axis. To see this, it is sufficient to construct
a hyperbolic exponential tract on which $F_{\lambda,m}$ has asymptotic direction $\theta^*$. Let $B_r$ an open disk
containing $F_{\lambda,m}(-m)$, the preimage of this set is an open set similar to $H_{|\lambda|,m}$. Let $D$ the complement
of this set. We have that $F_{\lambda,m}$ maps $D$ onto the exterior of $B_r$, then $D$ is an exponential tract
for $F_{\lambda,m}$. We may choose the negative real axis to define the fundamental domains in $D$. Since the
curves $\sigma_k$ for $k \in  \bz$  are mapped by $F_{\lambda,m}$ onto this axis, it follows that $D$ has asymptotic
direction $\theta^* = 0$. Furthermore, since $F_{\lambda,m}(z) \, = \, \lambda z^m \exp(z)$, one may check readily that
$D$ is a hyperbolic exponential tract.
\end{proof}

\begin{figure}[hbt]
\psfrag{a}[][]{\small $A^*(0)$}
\psfrag{b}[][]{\small $W$}
\psfrag{c}[][]{\small $D$}
\psfrag{d}[][]{\tiny $m \, - \,1 $}
\psfrag{de}[][]{\small $\bd$}
\psfrag{fi}[][]{\small $\phi$}
\psfrag{fu}[][]{\small $F_{\la,m}$}
\psfrag{pol}[][]{\small $z \to z^m$}
\centerline{\includegraphics[width=0.7 \textwidth]{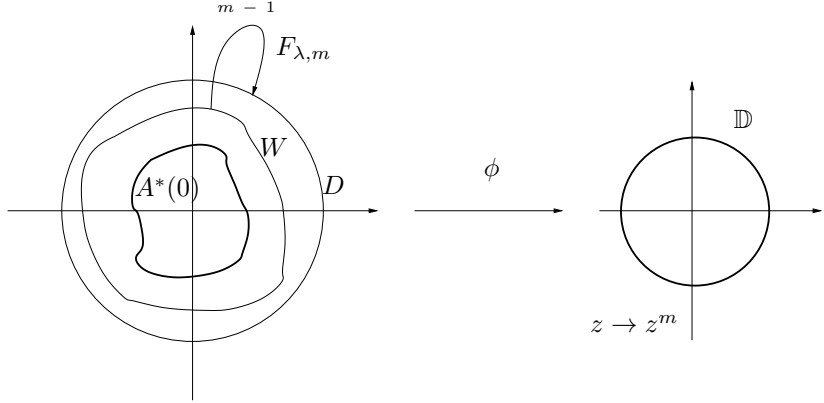}}
\caption{\small{ $F_{\lambda,m}$ is a polynomial-like mapping of degree $m$ near the origin.}}
\label{quasi}
\end{figure}

\begin{proposition}
If $\lambda \notin C_m^0$ then $A(0)$ has infinitely many components. Moreover, if $|\lambda| > (\frac{e}{m-1})^{m-1} $,  the boundary of $A^*(0)$ is a quasi-circle.
\end{proposition}
\begin{proof}

Using lemma \ref{pa00} we have that $A(0)$ has either one or infinitely many connected components. If we suppose
that $A(0)$ has only one connected component, then $A(0)$ is a completely invariant component of the Fatou set. We have
that all the critical values of $F_{\lambda,m}$ are in $A(0)$ (see \cite{Ba2}), and hence we conclude
that  $-m$ belongs to $A(0)$. However, is imposible if $\lambda \notin C_m^0$.

Let $\lambda \notin C_m^0$ such that $|\lambda| > (\frac{e}{m-1})^{m-1} $. The main idea of this proof is the same as that
used by Bergweiler in (\cite{B}). We will show that $F_{\la,m}$ is a polynomial-like of degree $m$ in a
neighbourhood of $0$, which includes the whole immediate basin. From the proof of proposition \ref{prim}, the disc
$D$ centered at $0$ of radius $m-1$ satisfies $\overline{D} \subset F_{\lambda,m}(D)$,
and hence $A^*(0) \subset D$.

Let $W$ be the component of $F^{-1}_{\lambda,m}(D)$ that contains the origin. It is clear that $\overline{W} \subset D$ and
$\overline{A^*(0)} \subset W$. Moreover, $F_{\lambda,m}$ is a proper function of degree $m$ from $W$ onto $D$, (see Fig. \ref{quasi}).
In the terminology of polynomial-like mappings, developed by Douady and Hubbard (\cite{DH}),
the triple $(F_{\lambda,m};W,D)$ is a polynomial-like mapping of degree $m$. By the Straightening theorem,
there exists a quasiconformal mapping, $\phi$, that conjugates  $F_{\lambda,m}$ to a polynomial $P$
of degree $m$, on the set $W$. That is $(\phi^{-1} \circ F_{\lambda,m} \circ \phi )(z)  = P(z) $ for all $z \in W$.
Since $z=0$ is superattracting for $ F_{\lambda,m}$ and $\phi$ is a conjugacy, we have that $z=0$ is
superattracting for $P$. Hence, after perhaps a holomorphic change of variables, we may assume that $P(z)=z^m$.

Hence, $\partial A^*(0) = \phi (\bp)$, and the theorem follows.
\end{proof}

\begin{remark}
The reason to ask for $|\la| > (\frac{e}{m-1})^{m-1}$ as a condition is as follows. We want to find
a value $K>0$ such that if $|z|=K$ then $|F_{\lambda,m}(z)| > K$. This condition is equivalent to

\[
|F_{\lambda,m}(z)| \ge  |\lambda| |z|^m e^{-|z|} = |\lambda| (K)^m e^{-K} > K
\]
\noindent or equivalently
\[
|\lambda| > K^{1-m}e^K.
\]
We want to use this argument for the largest possible region of values of $\lambda$. Hence,
we choose $K>0$, such that  $K^{1-m}e^K$ is minimum. This minimum value is reached exactly at $K=m-1$.
\end{remark}

\subsection{Proof of Theorem C}

In this Sec. we describe the capture zones $C_m^2$ (Proposition \ref{par_real} and
Proposition \ref{thi}) and $C_m^3$ (Proposition \ref{banpp}).

We will construct the open set $C_m^2$ in two steps. In the first (Proposition  \ref{par_real}) we obtain an
unbounded interval of real numbers $I$, such that for all  $\lambda \in I, \, F^2_{\lambda,m}(-m) \,$
lies in $\, A^*(0)$. In the second  (Proposition \ref{thi}), we will extend this construction to $\lambda \, $ in
$\, {\mathbb C}$. We denote $\lambda = \lambda_1 + \lambda_2 i$, where $\lambda_1 \hbox{ and }  \lambda_2$, are
the real and imaginary parts of $\lambda$.

\begin{proposition}\label{par_real}
There exists an unbounded interval,$\,I$, such that for all real numbers $\lambda_1 \in I$, we have that
$F_{\lambda_1,m}^2(-m) \in D_{\epsilon_0} \subset A^*(0)$.
\end{proposition}

\begin{proof}
Hereafter we denote  $r_m = (\frac{e}{m})^m$. We take  $ \lambda_1  \in  {\mathbb R} $, and we impose that
$F_{\lambda_1,m}(-m) \in H_{|\lambda_1|,m}$.

If we calculate  $F_{\lambda_1,m}(-m)$ we obtain
\[
F_{\lambda_1,m}(-m) \, = \, \lambda_1 (-m)^m \exp(-m) \, = \, (-1)^m \frac{\lambda_1}{r_m}.
\]

\noindent This real value lies in $H_{|\lambda_1|,m}$, if and only if

\[
|h(F_{\lambda_1,m}(-m))| \, < \, \frac{\epsilon_0}{|\lambda_1|}.
\]

\noindent Recall that  $h(x) = x^m \exp(x)$, and $\epsilon_0$ only depends on $|\lambda_1|$ and $m$. This condition
is equivalent to

\[
\frac{|\lambda_1|^{m+1}}{r_m^m} \exp((-1)^m \frac{\lambda_1}{r_m}) \, < \, \epsilon_0.
\]

Using lemma \ref{lem_e0} we have that $\epsilon_0 \, \ge \, \min\{1,(\frac{1}{|\lambda_1|e})^{1/(m-1)}\}$. If we use this
explicit lower bound we may impose

\[
\frac{|\lambda_1|^{m+1}}{r_m^m} \exp((-1)^m \frac{\lambda_1}{r_m}) \, < \,\min\{1,(\frac{1}{|\lambda_1|e})^{1/(m-1)}\}
\]

We define the auxiliary function

\[
l(\lambda_1) \, = \,  \left\{
\begin{array}{ll}
|\lambda_1|^{m+1+\frac{1}{m-1}} \, \exp((-1)^m \, \frac{\lambda_1}{r_m}) \exp(1/(m-1))  & \,  \hbox{if } \, |\lambda_1| > 1/e   \\
|\lambda_1|^{m+1} \exp((-1)^m \, \frac{\lambda_1}{r_m}) & \, \hbox{if } \, |\lambda_1| \le 1/e
  \end{array}\right.
\]

\noindent and the above inequality is transformed into $l(\lambda_1)< r_m^m$.

Using some elementary properties of function $l(\lambda_1)$, one can see that
\[
\left\{
\begin{array}{ll}
lim_{\lambda_1 \to -\infty} l(\lambda_1) \, = \, 0 & \, \hbox{if } $m$ \hbox { is even}\\
lim_{\lambda_1 \to +\infty} l(\lambda_1) \, = \, 0 & \, \hbox{if } $m$  \hbox { is odd}.
\end{array}\right.
\]
Since $l(\lambda_1)$ is continous and positive and it has a finite number of relative maxima and minima,
we can find an unbounded interval of  real numbers such that $ l(\lambda_1)< r_m^m $.

If $m$ is even, we define  $I=(-\infty,-D_0)$, where  $-D_0=-D_0(m) \le 0$, is the smallest of the values such that $l(\lambda_1)=r_m^m$. If $m$
is odd, we choose  $I=(D_0,+\infty)$, with $D_0 =D_0(m)\ge 0$, such that $D_0$ is the largest of the values
for which $l(\lambda_1)=r_m^m$.
\end{proof}

\begin{proposition}\label{thi}
Let $D_0(m) > 0$ be as in Proposition \ref{par_real}.There exists a function  $\alpha=\alpha(|\lambda|,m) \in 
(\pi/2,\pi)$, such that
\[
\begin{array}{ll}
\bullet \hbox{ for $m$ even,  the set $C_m^2$ contains the open set} &  \left\{ \lambda \in {\mathbb C} \, \left|
\begin{array}{ll}
 |\lambda| & \, >  \,  D_0   \\
|Arg(\lambda)|  & \, > \, \alpha    \end{array}\right. \right\} \\
\bullet \hbox{ for $m$ odd,  the set $C_m^2$ contains the open set} & \left\{ \lambda \in {\mathbb C} \, \left|
\begin{array}{ll}
|\lambda|  & \, > \, D_0   \\
 |Arg(\lambda)| & \, < \,  \pi - \alpha  \end{array}\right. \right\}
\end{array}
\]
\end{proposition}

\begin{figure}[hbt]
\label{}
\psfrag{x}[][]{ $\la_1$}
\psfrag{y}[][]{$\la_2$}
\psfrag{s}[][]{ $S$}
\psfrag{s1}[][]{$\widetilde{S}$}
\psfrag{h}[][]{ $H_S$}
\psfrag{a}[][]{ $\alpha$}
\psfrag{z}[][]{ $\zeta$}
\psfrag{i}[][]{ $\la_1^*$}
\psfrag{z1}[][]{ $\bar{\zeta}$}
\centerline{\includegraphics[width=0.5 \textwidth]{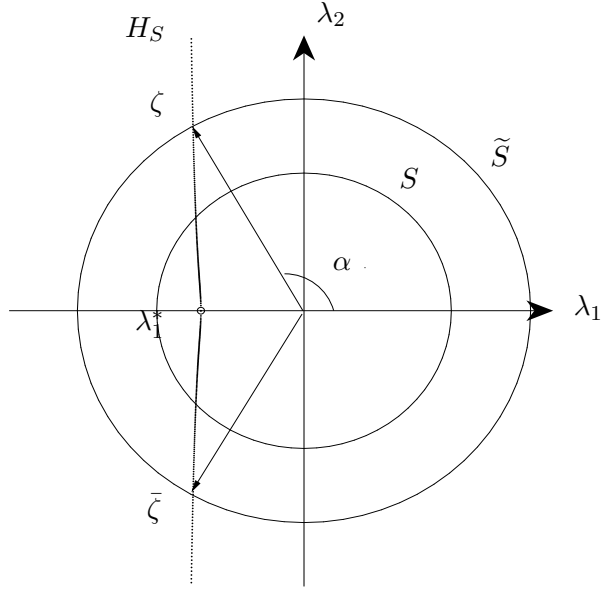}}
\caption{\small{Sketch of construction of the set $H_m, \, m=2.$ }}
\label{para}
\end{figure}

\begin{proof}
Given $\lambda_1^*  \in I$, we denote by $S$ the circle of radius $|\lambda_1^*|$ and centered at the origin. We will 
find all complex
numbers   $\lambda \, $ in $\, S$, such that  $F_{\lambda,m}(-m) \in H_{|\lambda|,m}$.

All complex numbers  $\la \in S$ have the same  $H_{|\lambda|,m}$ set, since this set only depends on $|\lambda|$ and $m$. We denote it by $H_S$.

When $\lambda $ belongs to $S$, the image of the critical point, $F_{\lambda,m}(-m)=(-1)^m \frac{\lambda}{r_m}$, 
belongs to another circle, namely $\widetilde{S}$, and its argument verifies 
\[
Arg(F_{\lambda,m}(-m)) \, = \, 
\left\{
\begin{array}{ll}
Arg(\la)  & \,  \hbox{ if } m \hbox{ is even}   \\
Arg(\la)+ \pi  & \,  \hbox{ if } m \hbox{ is odd}   \\
  \end{array}\right.
\]

This circle is concentric with respect to $S$ and its radius is
equal to $\frac{|\lambda_1^*|}{r_m}$ (see Fig. \ref{para}), which is larger
than the radius of $S$ if $m=2$, and smaller if $m \ge 3$.

We take $\lambda_1^* \,$ on $\,  I$. Using the construction above of the interval I, we obtain that $F_{\lambda_1^*,m}(-m) \in H_{|\lambda_1^*|,m}
=H_S$. This fact assures a non-empty intersection of  $\partial H_S$ with $\widetilde{S}$.

Using the analytic definition of $H_S$ (proof of Proposition \ref{th}), we can calculate  $\partial H_S \cap  \widetilde{S}$.
We find this intersection by solving:

\[
 \left\{
\begin{array}{ll}
\sqrt{ \la_1^2 \, + \, \la_2^2} \,  = \, \frac{|\lambda_1^*|}{r_m}  & \,  \la_1+i \la_2  \in  \widetilde{S}   \\
 \la_2 = +C(\la_1) =  \sqrt{ \left[ \frac{\epsilon_0}
{|\la_1^*|} \right]^{2/m} \exp(-\la_1 \frac{2}{m}) -  \la_1^2 } & \, \la_1 +i \la_2  \in \partial H_S
  \end{array}\right.
\]

It is not difficult to show that this system has two conjugate solutions namely $\zeta$ and $\bar{\zeta}$. If
we write $\zeta=\widetilde{\la_1} + i \widetilde{\la_2}$, then
\[
\widetilde{\la_1} = \ln{\frac{\epsilon_0 r_m^m}{|\lambda_1^*|^{m+1}}} \qquad  \widetilde{\la_2} = +C(\widetilde{\la_1})
\]

Let $\alpha=\alpha(|\la_1^*|,m) \in (\pi/2,\pi)\,$ be the argument of $\zeta$ (see Fig. \ref{para}).

If $m$ is even then for all complex numbers with modulus equal to $|\la_1^*|$, where $ \la_1^* \,$
lies in $\, I$, and argument greater than $\alpha$ in absolute value, it is verified that $F_{\lambda,m}(-m)
\in H_{|\lambda|,m}$. If $m$ is odd, the same is true for all complex numbers with modulus equal
to $|\la_1^*|$, where $\la_1^* \, $ lies in  $\, I$, and argument, in absolute value, smaller than $\pi-\alpha$.
\end{proof}

Parallel to the construction in dynamical plane we will now show the existence of a countable
number of horizontal bands, all of which are also capture zones. See Fig. \ref{pp2}- \ref{ppm}. Apparently, when $m$ is even, these
strips extend to $+\infty$, while if $m$ is odd they extend  to
$-\infty$. It also seems that  their width decreases as $m$ increases.

In the Sec. above we constructed similar strips around the curves $\sigma_k$. These curves were preimages of the negative
real axis. In this case, we define
\[
\Gamma_k = \{ \la \in \bc \hbox{ such that }  F_{\la,m}(-m) \in  \sigma_k \}.
\]

To find an expression for the curves $\Gamma_k$ we first write

\[
 F_{\la,m}(-m) = \la (-m)^m \exp(-m) = (-1)^m \frac{\la}{r_m}
\]
where $r_m = (\frac{e}{m})^m$. Recall that $z \in \sigma_k$ when
\[
Arg(\lambda) + m Arg(z) + Im(z) = (2k+1)\pi.
\]

Hence we need that
\[
Arg(\la) + m Arg(F_{\la,m}(-m)) + Im(F_{\la,m}(-m)) = (2k+1)\pi.
\]
If $m$ is even then  $Arg(F_{\la,m}(-m)) = Arg(\la)$, while if $m$ is odd then $Arg(F_{\la,m}(-m)) = Arg(\la) + \pi$;
thus we obtain the condition for $F_{\la,m}(-m) \in \sigma_k$

\[
 \left\{
\begin{array}{ll}
Arg(\la) + m Arg(\la) + \frac{|\la|}{r_m} \sin(Arg(\la)) = (2k+1)\pi &\, \hbox{if } $m$ \hbox { is even}\\
Arg(\la) + m (Arg(\la)+\pi) + \frac{|\la|}{r_m} (-1)\sin(Arg(\la)) = (2k+1)\pi & \, \hbox{if } $m$  \hbox { is odd}
\end{array}\right.
\]

Solving for  $|\la|$, we obtain a function of $Arg(\la)$, which we denote by $\phi$. Explicitly, the curve $\Gamma_k$ can
be written as

\[
 \left\{
\begin{array}{ll}
|\la| = \phi(Arg(\la)) = r_m \frac{(2k+1)\pi - (m+1) Arg(\la)} {\sin(Arg(\la))}  \qquad -\pi \le Arg(\la) \le \pi & \hbox{if } $m$ \hbox{ is even} \\
 |\la| = \phi(Arg(\la)) = r_m \frac{(2k+1-m)\pi - (m+1) Arg(\la)} {-\sin(Arg(\la))}   \qquad -\pi \le Arg(\la) \le \pi & \hbox{if } $m$ \hbox{ is odd}
\end{array}\right.
\]
As in the Sec. above, we need to impose $\phi(Arg(\la)) \ge 0$. If we denote $\theta = Arg(\la)$, we have:

\noindent if $m=2j \hbox{ for } j \in \bz$ 
\[
\Gamma_k = \phi(\theta) e^{i\theta}   \left\{  \begin{array}{llll}   0< \theta < \pi \qquad & \hbox{  if  }\quad
 k \ge j+1 \\
0< \theta <\frac{(2k+1)\pi}{m+1} \qquad & \hbox{  if  }\quad 0 \le  k \le j \\
 \frac{(2k+1)\pi}{m+1}< \theta <0 \qquad & \hbox{  if  }\quad -(j+1) \le  k \le 0 \\
 -\pi< \theta <0  \qquad & \hbox{  if  }\quad k \le -(j+2) \\ \end{array}\right. ;
\]

\noindent if  $m=2j+1 \hbox { for } j \in \bz$ 
\[
\Gamma_k =  \phi(\theta) e^{i\theta}   \left\{  \begin{array}{lllll}  0< \theta < \pi \qquad & \hbox{  si  }\quad  k \ge m \\
-\pi < \theta < 0  \, \cup \,  \frac{2k+1-m}{m+1} < \theta < \pi   \qquad & \hbox{  if  }\quad j+1 \le  k \le m-1 \\
-\pi< \theta < \pi  \qquad & \hbox{  if  }\quad   k = j \\
 -\pi< \theta < \frac{2k+1-m}{m+1} \, \cup \, 0 < \theta < \pi  \qquad & \hbox{  if  }\quad j-1 \le  k \le 0 \\
 -\pi< \theta <0  \qquad & \hbox{  if  }\quad k \le -1 \\ \end{array}\right.
\]

\begin{figure}[hbt]
	\centering
	\subfigure[\scriptsize{Graph of $\Gamma_k$ for $m=3$} ]{
\psfrag{a}[][]{\scriptsize $\Gamma_4$}
\psfrag{b}[][]{\scriptsize $\Gamma_3$}
\psfrag{c}[][]{\scriptsize $\Gamma_2$}
\psfrag{d}[][]{\scriptsize $\Gamma_0$}
\psfrag{e}[][]{\scriptsize $\Gamma_{-1}$}
\psfrag{m}[][]{\scriptsize $\Gamma_{-}$}
\psfrag{f}[][]{\scriptsize $\Gamma_{-2}$}
\psfrag{g}[][]{\scriptsize $\Gamma_{-3}$}
\psfrag{h}[][]{\scriptsize $\Gamma_1$}
\psfrag{i}[][]{\scriptsize $\Gamma_0$}
\psfrag{j}[][]{\scriptsize $\Gamma_2$}
\psfrag{k}[][]{\scriptsize $\Gamma_1$}
	\fbox{ \includegraphics[width=0.3 \textwidth]{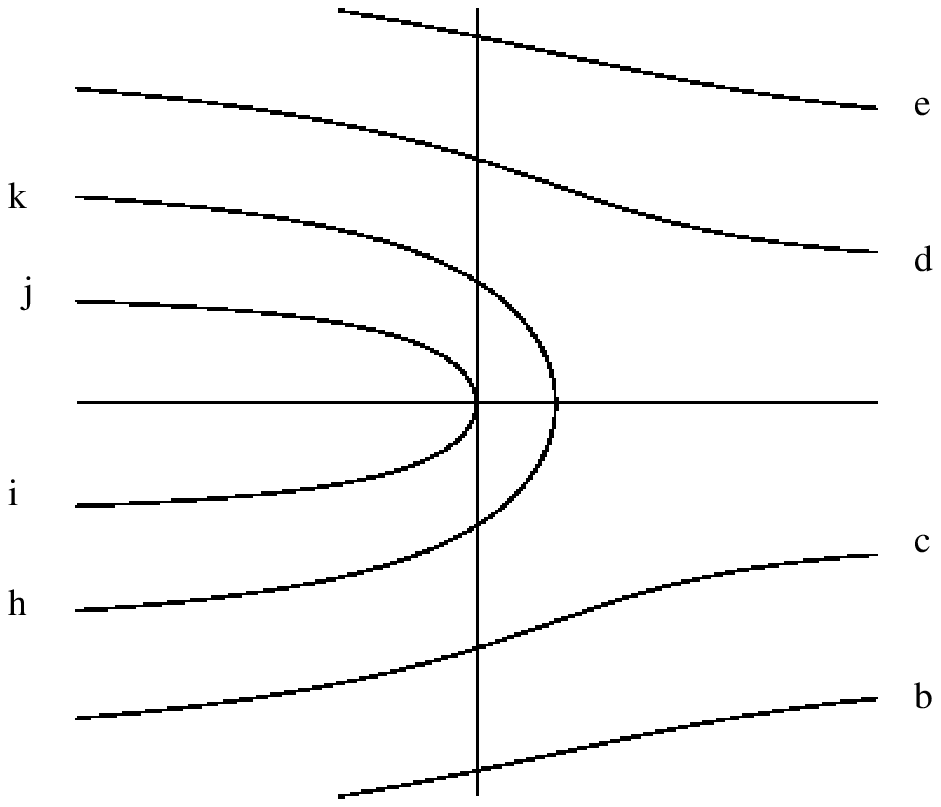}}}
	\hspace{0.5in}
	\subfigure[\scriptsize{The parameter plane of $F_{\la,3}$}  ]{
	\includegraphics[width=0.3 \textwidth]{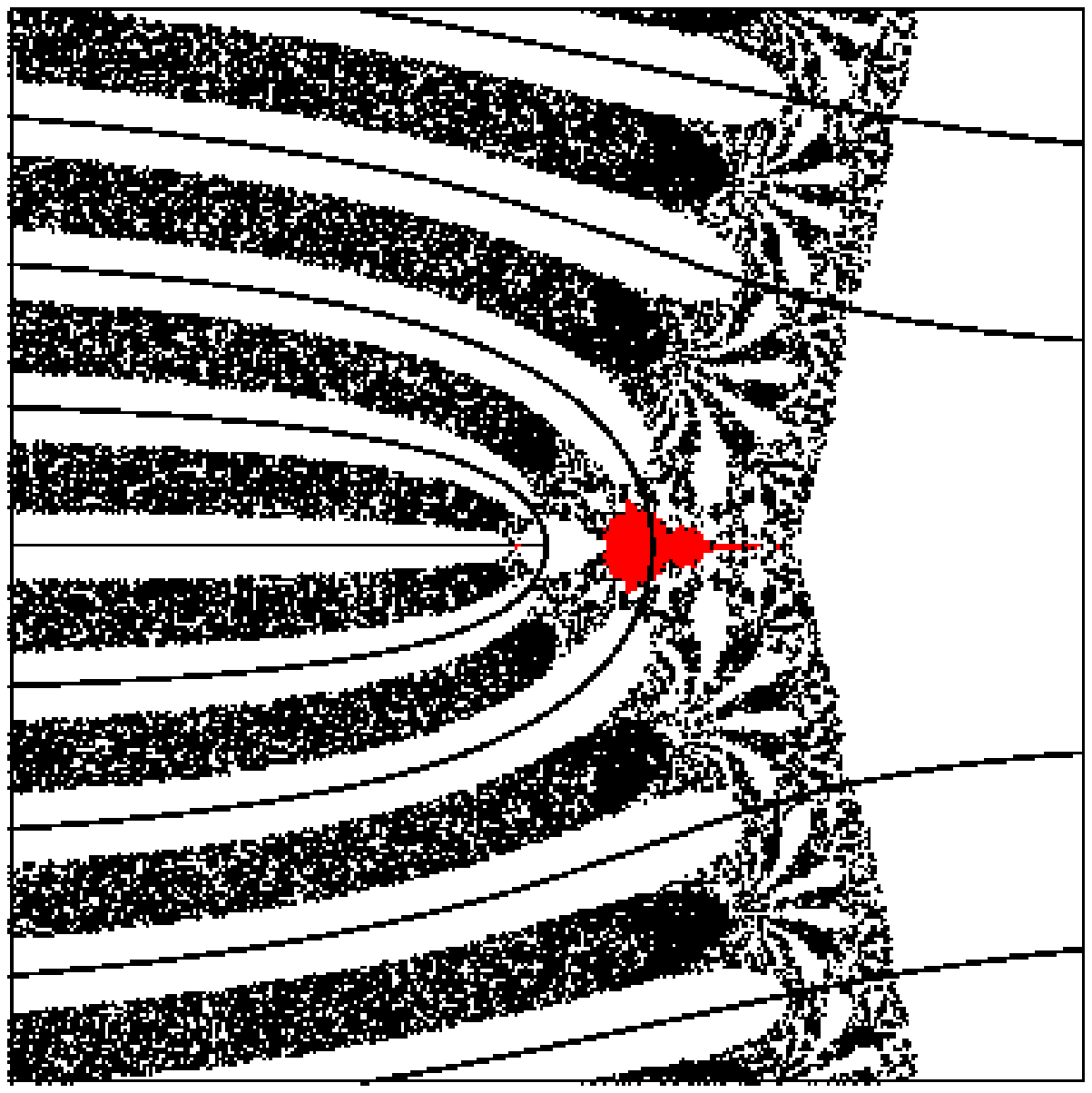}}
	\hspace{0.5in}
	\subfigure[\scriptsize{ Graph of $\Gamma_k$ for $m=4$}  ]{
\psfrag{a}[][]{\scriptsize $\Gamma_{-3}$}
\psfrag{b}[][]{\scriptsize $\Gamma_{-2}$}
\psfrag{c}[][]{\scriptsize $\Gamma_{-1}$}
\psfrag{d}[][]{\scriptsize $\Gamma_0$}
\psfrag{e}[][]{\scriptsize $\Gamma_1$}
\psfrag{f}[][]{\scriptsize $\Gamma_2$}
\psfrag{p}[][]{\scriptsize $\Gamma_{-6}$}
\psfrag{g}[][]{\scriptsize $\Gamma_{-5}$}
\psfrag{h}[][]{\scriptsize $\Gamma_{-4}$}
\psfrag{i}[][]{\scriptsize $\Gamma_3$}
\psfrag{j}[][]{\scriptsize $\Gamma_4$}
\psfrag{k}[][]{\scriptsize $\Gamma_5$}
	\fbox{ \includegraphics[width=0.3 \textwidth]{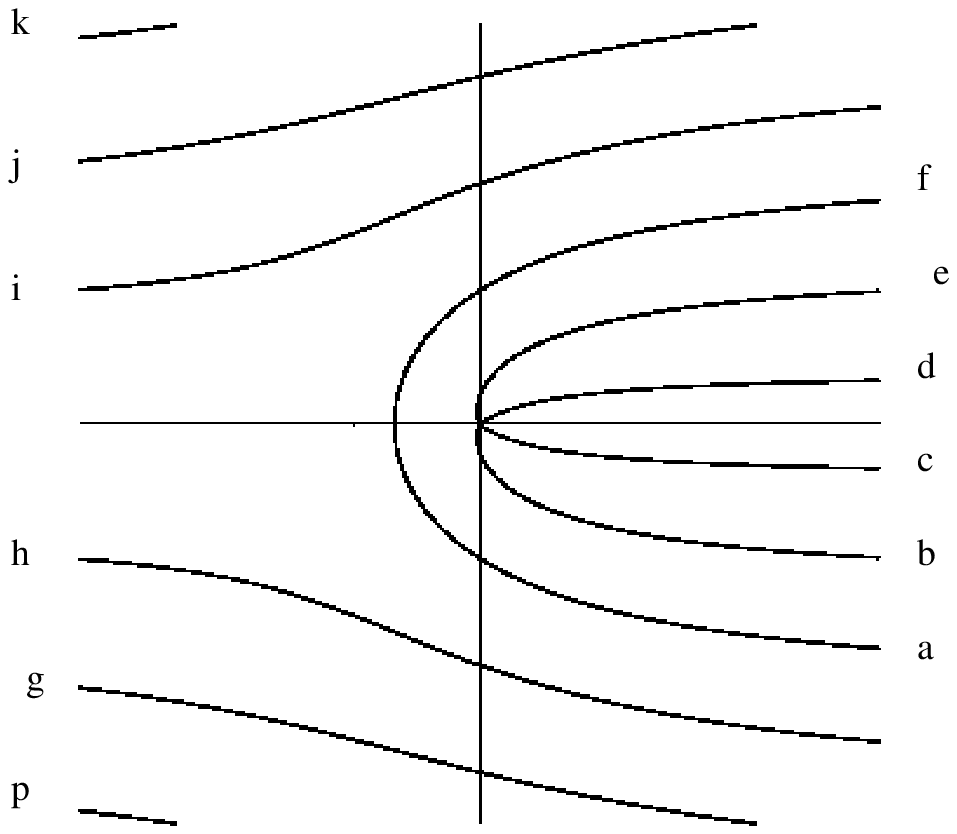}}}
	\hspace{0.5in}
	\subfigure[\scriptsize{Parameter plane of $F_{\la,4}$}  ]{
	 \includegraphics[width=0.3 \textwidth]{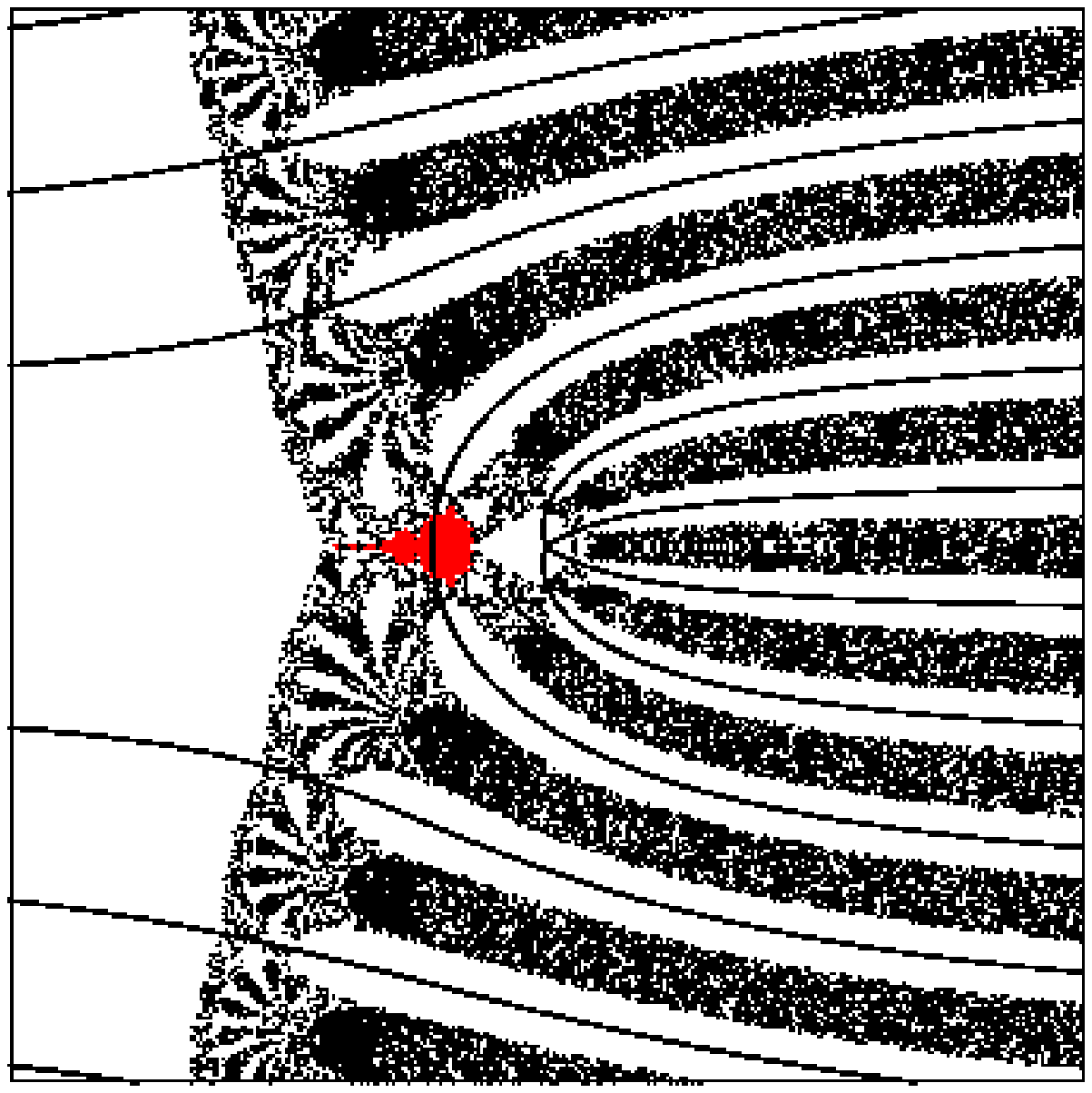}}
	\caption{\small{Strips in the parameter plane.}}
	\label{gamma}
\end{figure}

In Fig. \ref{gamma} we show some of these curves for some values of $m$. If we suppose that $m=2j$ is even, then
each $\Gamma_k$ tends asympotically to the line $Im(z) = (2k+1)\pi \, r_m$ as its real part tends to $+\infty$. We can
classify these curves in three types. The first one is formed by curves whose real part runs from $-\infty$ to $+\infty$. There
are two curves of the second kind, $\Gamma_j \hbox { and } \Gamma_{-(j+1)}$, with real part in $[-(m+1)r_m,\infty]$. The
third group is formed by $m$ curves, starting at the origin and tending to $+\infty$. These $m$ curves
have indexes between $j-1 \hbox { and } -j$.

If we take $m$ an odd index ($m=2j+1$), these curves tend asymptotically to the lines $Im(z) = 2(k-m)\pi \, r_m$ as their real part
tend to $-\infty$. As above, we can classify these curves in three types. The first one is formed by curves that extend from $-\infty$ to $+\infty$. The
second one is formed by the curve $\Gamma_j $, has a horseshoe shape, and cuts the real axis at the point
$(m+1)r_m$. The third one is formed by $m-1$ curves, starting at the origin and tending to $-\infty$. These $m$ curves
have indexes between $0 \hbox { and } m-1$, except for $\Gamma_j$.

Hence we have obtained some curves $\Gamma_k$ such that if $\la \in \Gamma_k$, then
$F_{\la,m}(-m) \in \sigma_k $. Hence, choosing $\la \in \Gamma_k$ with $Re(F_{\la,m}(-m))$ large enough, we obtain that
$F_{\la,m}(F_{\la,m}(-m)) \in  H_{|\la|,m}$. Since $Re(F_{\la,m}(-m))= Re(\la)(-1)^m r_m$, this corresponds to taking $Re(\la)$ or $-Re(\la)$ large enough depending on $m$
being even or odd.

By construction, the half curves we just defined belong each to $C_m^3$. We will
now show that a neighbourhood of $\Gamma_k$ of asymptotic width equal to $r_m \pi$
is also part of $C_m^3$. We fix a value $k \in {\mathbb Z}$ and we suppose that
$\la = \la_1 + i\la_2$, where $\la_1 \hbox{ and } \la_2$ are real numbers. We will prove the following result.

\begin{proposition} \label{banpp}
If  $m$ is even, for all $ \la_2 \in (r_m (\frac{\pi}{2}+2k\pi)  \, ,
\, r_m(\frac{3\pi}{2}  + 2k\pi)) $
 there exists  $\la_1^* \,$ such that, for all $ \la_1 > \la_1^* \,$ then $ F_{\lambda,m}^3(-m) \in A^*(0)$.

If $m$ is odd, for all $ \la_2 \in  (r_m(-\frac{\pi}{2} + 2k\pi)  \, , \, r_m(\frac{\pi}{2}  +  2k\pi)) $
there exists $ \la_1^* \,$  such that, for all $  \la_1 < \la_1^* \,$ then $F_{\lambda,m}^3 (-m) \in A^*(0)$.

\end{proposition}

\begin{proof}
Assume that $m$ is even (the odd case is completely symetric), and let $\la_2 \in (r_m (\frac{\pi}{2}+2k\pi)  \, , \, r_m(\frac{3\pi}{2}  +
2k\pi)) $. We recall
that proposition \ref{thbandas} assures that, for all $\, y \in ( (2k+1) \pi - Arg( \lambda) - \frac{\pi}{2} \,  , \, (2k+1) \pi -
Arg( \lambda) + \frac{\pi}{2}) , \, $ there exists $\, x_* \in  {\mathbb R} $ such that, for all $\,x \ge x_* $, the point
$F_{\lambda,m} (x+yi) \in H_{|\lambda|,m}$.

Using that $F_{\lambda,m} (-m)\, = \, \frac{\la_1}{r_m}  \, + \,  \frac{\la_2}{r_m} i $, it suffices to prove that
\[
Im(F_{\la,m}(-m)) = \frac{\la_2}{r_m}  \in (\frac{\pi}{2} + 2k\pi - Arg(\lambda) \, , \, \frac{3\pi}{2}  + 2k\pi -  Arg(\lambda) ).
\]
\noindent Choosing
\[
Re(F_{\la,m}(-m)) = \frac{\la_1}{r_m} > x^*.
\]
\noindent The first condition is equivalent to $Arg(\lambda) \in (\alpha_1 , \alpha_2)$,

\begin{figure}[hbt]
\psfrag{a1}[][]{ \scriptsize{$r_m(\frac{3\pi}{2}+2k\pi) i$}}
\psfrag{a2}[][]{\scriptsize{$r_m(\frac{\pi}{2}+2k\pi) i$}}
\psfrag{yr}[][]{ \scriptsize{$\la_2 / r_m i$}}
\psfrag{j}[][]{$\alpha_2$}
\psfrag{l}[][]{ $\alpha_1$}
\psfrag{al}[][]{ $Arg(\lambda)$}
\psfrag{rp}[][]{ $\scriptsize{\widetilde{\la_1}}$}
\psfrag{c1}[][]{\small{$\alpha_2 > 0$}}
\psfrag{c2}[][]{$\alpha_1 < 0$}
\psfrag{r}[][]{$r$}
\psfrag{s}[][]{$s$}
\centerline{\includegraphics[width=0.4 \textwidth]{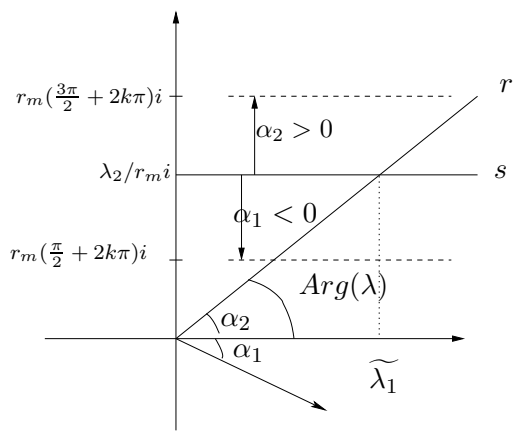}}
\caption{\small{Construction of the value  $\widetilde{\la_1}$  }}
\label{ban}
\end{figure}

\noindent where $\alpha_1 = \frac{\pi}{2}  + 2k\pi - \frac{\la_2}{r_m} < 0 \,\hbox{ and } \, \alpha_2 = \frac{3\pi}{2}  + 2k\pi - \frac{\la_2}{m} >0$ (Fig. \ref{ban}).

Suppose that $k>0$. We denote by $r$ the line through the origin with slope  $tan(\alpha_2)$, and let $s$ be the horizontal
line through $\frac{\la_2}{r_m} i$. We also denote by $\widetilde{\la_1}$ the abscissa of the intersection
point between the lines $r \, \hbox{ and } \, s $ (Fig. \ref{ban}). All values of $\la$ on  $s$, with abscissa greater than  $\widetilde{\la_1}$ verify that $0< arg(\la) < \alpha_2$. Finally,
we define $\la_1^* = max\{r_m x_* ,r_m \widetilde{\la_1}\}$, and for this value both conditions are verified.
Therefore, $F^2_{\la,m} (-m) \in H_{|\la|,m}$, and it follows that $F^3_{\la,m} (-m) \in D_{\epsilon_0}
\subset A^*(0)$. If $k \le 0$, we replace the line of slope $tan(\alpha_2)$, by a line with slope
$tan(\alpha_1)$.

\end{proof}

\noindent{\bf Acknowledgments}

We wish to thank Robert Devaney and Xavier Jarque for very helpful 
discussions.


\end{document}